\documentclass[twocolumn,journal]{IEEEtran}
\usepackage{hyperref}
\usepackage{siunitx}
\usepackage{makecell}
\usepackage{easyReview}
\usepackage{amssymb}
\usepackage{amsmath}
\usepackage{amsfonts}
\usepackage{graphicx}
\usepackage{color}
\usepackage{caption}
\usepackage{multirow}
\usepackage{diagbox}
\usepackage{lineno}
\usepackage{subcaption}
\usepackage{mathrsfs}
\usepackage{amscd}
\usepackage{url}
\usepackage{enumitem}
\usepackage{booktabs}
\usepackage{tabularx}
\usepackage{makecell}
\usepackage{bm}
\usepackage{bbm}
\graphicspath{{Fig/}}  
\usepackage{tikz}
\usetikzlibrary{shadows, positioning, fit, backgrounds, calc, shapes.geometric, arrows.meta}
\usepackage{dashrule}
\newtheorem{theorem}{Theorem}[section]

\newtheorem{algorithm}{Algorithm}[section]

\newtheorem{corollary}{Corollary}[section]

\newtheorem{defi}{Definition}[section]

\newtheorem{lem}{Lemma}[section]

\newtheorem{remark}{Remark}[section]

\newenvironment{proof}[1][Proof]{\noindent\textbf{#1.} }{\ \hfill\rule{0.3em}{0.5em}}

\usepackage{mathtools}

  \newcommand{\A}{\mathcal{A}}
\newcommand{\mB}{\mathcal{B}}  
\newcommand{\mC}{\mathcal{C}}  
     
\newcommand{\E}{\mathcal{E}}   \newcommand{\mE}{\mathcal{E}}
\newcommand{\F}{\mathcal{F}}
\newcommand{\G}{\mathcal{G}}

\newcommand{\I}{\mathcal{I}}

\newcommand{\K}{\mathcal{K}}
\newcommand{\mL}{\mathcal{L}}  
\newcommand{\mM}{\mathcal{M}}

\renewcommand{\P}{\mathcal{P}} 
\newcommand{\Q}{\mathcal{Q}}   
\newcommand{\R}{\mathbb{R}}    \newcommand{\mR}{\mathcal{R}}
\newcommand{\mS}{\mathcal{S}}
\newcommand{\mT}{\mathcal{T}}
\newcommand{\U}{\mathcal{U}}
\newcommand{\V}{\mathcal{V}}
  \newcommand{\W}{\mathcal{W}}
\newcommand{\X}{\mathcal{X}}
\newcommand{\Y}{\mathcal{Y}}
\newcommand{\Z}{\mathcal{Z}}

\newcommand{\rank}{\operatorname{rank}}

\newcommand{\TF}{\operatorname{TF}}
\newcommand{\HS}{\operatorname{HS}}
\newcommand{\RF}{\operatorname{RF}}
\newcommand{\range}{\operatorname{range}}

\usepackage{algorithm}
\usepackage{algpseudocode}
\usepackage[algo2e,ruled,linesnumbered]{algorithm2e}

\ifCLASSINFOpdf
\else
\fi
\begin{document}

\title{A Multi-View Coupled Tensor Decomposition for Lightweight Online Adaptive Traffic Prediction}

%
%
%
\author{
	Quan Yu,
	Jie Ni,
	Yu-Hong Dai,
	Xiongjun Zhang\thanks{*Corresponding author: Yu-Hong Dai.}%
	\thanks{This work was supported  in
		part by the National Key R\&D Program of China under Grant 2021YFA1000300, and Grant 2021YFA1000301; in part by the National Natural Science Foundation of China under Grant 12171189, Grant 12021001, and Grant 92473208; in part by the Hubei Provincial Natural Science Foundation of China under Grant
		2025AFB966; and in part by the Fundamental Research Funds for the Central
		Universities under Grant XJ2026000601.}
	\thanks{Quan Yu, and Xiongjun Zhang are with the School of Mathematics and Statistics, Central China Normal University, Wuhan 430079, China, and also with the Key Laboratory of Nonlinear Analysis and Applications (Ministry of Education), Central China Normal University (e-mails: quanyu@ccnu.edu.cn, xjzhang@ccnu.edu.cn).}
	\thanks{Jie Ni is with the School of Mathematics and Statistics, Beijing Institute of Technology, Beijing 100081, China (e-mail: jie\_ni2@163.com).}
	\thanks{Yu-Hong Dai is with the State Key Laboratory of Mathematical Sciences, Academy of Mathematics and Systems Science, Chinese Academy of Sciences, Beijing 100190, China, and also with the School of Mathematical Sciences, University of Chinese Academy of Sciences, Beijing 100049, China (e-mail: dyh@lsec.cc.ac.cn).}%
}

%
%

\markboth{}
{Shell \MakeLowercase{\textit{et al.}}: Bare Demo of IEEEtran.cls for IEEE Journals}
%



\maketitle

\begin{abstract}
	
Accurate online traffic prediction is essential for intelligent transportation systems, where forecasting must be performed continuously under imperfect sensing conditions. Missing observations and anomalous disturbances make this task challenging, particularly when prediction relies on a single traffic view. This paper proposes a Multi-View Coupled Tensor Decomposition (MVCTD) model for online traffic prediction from imperfect multi-view observations, such as speed, flow, and occupancy. The proposed model uses coupled tensor decomposition to build a structured latent forecasting space, in which shared spatial structures across traffic views and view-specific temporal dynamics are jointly modeled. A group sparse regularization is further introduced to capture correlated abnormal responses induced by real traffic anomalies and thus reduce their influence on forecasts. For streaming deployment, MVCTD performs iterative refinement only on the current latent tensor, while the remaining model variables are updated by lightweight closed-form steps based on summarized historical information, thereby avoiding repeated optimization over the full historical sequence. Experiments on real-world traffic datasets demonstrate that MVCTD achieves accurate forecasts with favorable runtime under severe missingness, confirming its suitability for online traffic prediction.

\end{abstract}


%
\IEEEpeerreviewmaketitle

\section{Introduction}
\IEEEPARstart{O}{nline} traffic prediction is a fundamental task in Intelligent Transportation Systems (ITS), where traffic states must be continuously inferred from streaming measurements and extrapolated over future horizons \cite{XLCL23}. Accurate forecasts support proactive traffic control, congestion mitigation, and route guidance \cite{AN17}. Modern sensing and communication infrastructures provide large-scale traffic measurements from urban road networks, but these data are high-dimensional, strongly spatio-temporal, and often available only in an imperfect form. As a result, effective online prediction requires not only modeling recurrent traffic dynamics, but also coping with missing observations, abnormal disturbances, and stringent computational constraints \cite{CS21}.

Traffic forecasting has been studied extensively in the literature. Classical time-series models, including Autoregressive (AR), Autoregressive Integrated Moving Average (ARIMA), and Seasonal ARIMA, model future values from historical temporal dependencies and have been applied to short-term traffic flow and occupancy prediction \cite{XLCL23,AC79,HAS95,LF99,WDB98}. However, these models are primarily designed for univariate or low-dimensional series and therefore have limited ability to capture network-wide spatial correlations. Vector Autoregression (VAR) extends autoregressive modeling to multiple correlated variables, and low-rank factorization has further improved its scalability for spatio-temporal data. For example, Bayesian Temporal Matrix Factorization (BTMF) combines matrix factorization with a Gaussian VAR process to regularize temporal factors and improve prediction under missing data \cite{CS21}. Related low-rank and dynamical system models have also been used in broader spatio-temporal forecasting applications \cite{TKU17,JSJZ19,WZ26,LSWZ25}.

Deep learning methods provide another major line of work by learning nonlinear temporal and spatial dependencies from data \cite{YWW22}. Representative examples include deep belief networks \cite{HSHX14}, recurrent neural networks and their gated variants such as long short-term memory and gated recurrent unit models \cite{Yas99,MMR16,MTW15,SVL14,FWL22,YPL21,FZL16}, convolutional neural networks \cite{ZZQ18,LLJ19,DJC19}, and graph convolutional networks that explicitly incorporate road-network topology \cite{AUA25,CC22,GLF19,HYYX23,WZCY22,ZHD22}. These models are powerful when sufficient clean training data and computational resources are available. Nevertheless, their performance can deteriorate when online inputs are heavily incomplete or contaminated, and the cost of updating the model for newly arriving data may be high in streaming deployment.

Despite these advances, online prediction from real-world traffic streams remains challenging for three reasons. First, traffic observations are inherently multi-view: variables such as speed, flow, and occupancy describe the same evolving network state from complementary perspectives \cite{HJJ24,YKB14}. Forecasting each view separately may ignore useful cross-view information, whereas naively stacking all variables may obscure view-specific temporal behavior. Second, streaming observations are often incomplete and occasionally corrupted by accidents, sensor faults, or other irregular events \cite{FZWZ22,LZC24}. Treating imputation, anomaly filtering, and forecasting as separate preprocessing and prediction stages can propagate estimation errors into the final forecast. Third, online deployment requires fast adaptation after new partial observations arrive. Repeatedly re-optimizing the entire historical sequence is computationally unattractive, especially when the history length grows over time.

These issues motivate a forecasting framework that couples multiple traffic views, learns from incomplete observations, separates regular dynamics from anomalous disturbances, and updates efficiently in a streaming setting. To this end, we propose a Multi-View Coupled Tensor Decomposition (MVCTD) model for lightweight online adaptive traffic prediction. MVCTD represents daily traffic states as third-order tensors indexed by road segment, time interval, and traffic view, and decomposes the regular traffic component into a shared spatial dictionary and day-specific latent temporal tensors through coupled tensor factorization. Forecasting is then performed in the latent space, where shared cross-view structures and view-dependent temporal dynamics are modeled together with temporal regularity and sparse abnormal components. For streaming deployment, MVCTD adopts a two-stage scheme that initializes the model from historical incomplete data and then performs lightweight online adaptation by refining only the current latent tensor while updating the remaining variables through closed-form steps based on summarized historical information.

The main contributions of this paper are summarized as follows:
\begin{itemize}
	\item We propose MVCTD, a multi-view coupled tensor decomposition framework that constructs a structured latent forecasting space for online traffic prediction. It couples traffic views through a shared spatial dictionary and view-dependent latent temporal factors, enabling cross-view complementarities to be exploited while retaining view-specific temporal dynamics.

	\item We develop a robust latent learning formulation that recovers predictive traffic states from incomplete and anomaly corrupted observations within a single optimization problem. The formulation combines latent space autoregression with cross-view group sparse anomaly modeling and latent temporal regularization for global periodicity and local smoothness, enabling recurrent traffic dynamics to be learned while reducing the risk of treating corrupted observations as regular traffic patterns.
	
	\item We design a lightweight two-stage optimization scheme for streaming deployment. In the offline stage, the model is initialized from historical incomplete data. In the online stage, only the current latent tensor is refined, while the remaining model variables are updated via closed-form solutions. This avoids repeated optimization over the full historical sequence, thereby reducing computational overhead and enabling efficient online adaptation.
		
\end{itemize}

The remaining parts of this paper are organized as follows. Section \ref{Sec:pre} introduces the tensor algebra and notation used throughout the paper. Section \ref{Sec:MVCTD} presents the proposed MVCTD model. Section \ref{Sec:Alg} develops the corresponding optimization algorithm. Section \ref{Sec:exp} reports the experimental results, and Section \ref{Sec:con} concludes the paper.
The proof of the main theorem is deferred to the supplementary material.

\section{Preliminaries}\label{Sec:pre}

This section introduces the notation and tensor operations used in the proposed MVCTD framework. We focus on the algebraic tools needed for coupled tensor factorization and latent temporal regularization.

For a positive integer $n$, let $[n]=\{1,2,\ldots,n\}$. Scalars, vectors, matrices, and tensors are denoted by lowercase letters, boldface lowercase letters, uppercase letters, and calligraphic letters, respectively. The fields of real and complex numbers are denoted by $\mathbb{R}$ and $\mathbb{C}$. For tensors $\mathcal{X},\mathcal{Y}\in\mathbb{R}^{n_1\times n_2\times n_3}$, let $\X_{ijk}$ and $\Y_{ijk}$ denote their $(i,j,k)$-th entries, respectively; the inner product is defined as $\langle \mathcal{X},\mathcal{Y}\rangle=\sum_{i=1}^{n_1}\sum_{j=1}^{n_2}\sum_{k=1}^{n_3}\X_{ijk}\Y_{ijk}$, with induced Frobenius norm $\|\mathcal{X}\|=\sqrt{\langle \mathcal{X},\mathcal{X}\rangle}$. For a matrix $X$, $\|X\|_*$ denotes its nuclear norm, i.e., the sum of its singular values. The Discrete Fourier Transform (DFT) of a vector $\bm{x}$ is denoted by $\hat{\bm{x}}$; for a tensor $\mathcal{X}$, $\hat{\mathcal{X}}$ denotes the tensor obtained by applying the DFT along the third dimension. The indicator function $\mathbb{I}_{\mathbb{S}}(x)$ equals one if $x\in\mathbb{S}$ and zero otherwise. The sign function is defined as $\operatorname{sgn}(x)=1$ if $x>0$, $\operatorname{sgn}(x)=0$ if $x=0$, and $\operatorname{sgn}(x)=-1$ if $x<0$. For $\bm{x},\bm{y}\in\mathbb{R}^n$, the circular convolution of $\bm{x}$ and $\bm{y}$ is defined as
$\bm{x}\bullet\bm{y}=\mathscr{C}(\bm{x})\bm{y}$, where $\mathscr{C}(\bm{x})$ is the circulant matrix given by
\begin{equation*}
	\mathscr{C}(\bm{x})=
	\begin{bmatrix}
		x_1 & x_n & \cdots & x_2 \\
		x_2 & x_1 & \cdots & x_3 \\
		\vdots & \vdots & \ddots & \vdots \\
		x_n & x_{n-1} & \cdots & x_1
	\end{bmatrix}.
\end{equation*}

\begin{defi}[t-product \cite{KM11}]\label{def:t-product}
	Let $\mathcal{X}\in\mathbb{R}^{n_1\times n_2\times n_3}$ and $\mathcal{Y}\in\mathbb{R}^{n_2\times n_4\times n_3}$. Their t-product $\mathcal{Z}=\mathcal{X}*\mathcal{Y}\in\mathbb{R}^{n_1\times n_4\times n_3}$ is defined by
	$\mathcal{Z}(i,j,:)=\sum_{k=1}^{n_2}\mathcal{X}(i,k,:)\bullet\mathcal{Y}(k,j,:)$.
\end{defi}
\begin{defi}[conjugate transpose \cite{KM11}]\label{def:ct}
	The conjugate transpose of $\mathcal{X}\in\mathbb{C}^{n_1\times n_2\times n_3}$, denoted by $\mathcal{X}^{\top}\in\mathbb{C}^{n_2\times n_1\times n_3}$, is obtained by conjugate transposing each frontal slice and reversing the order of the frontal slices from the second to the last.
\end{defi}
\begin{defi}[f-diagonal tensor \cite{KM11}]
	A tensor is f-diagonal if each of its frontal slices is a diagonal matrix.
\end{defi}
\begin{defi}[identity tensor \cite{KM11}]\label{def:it}
	The identity tensor $\mathcal{I}\in\mathbb{R}^{n\times n\times n_3}$ has the identity matrix as its first frontal slice and zero matrices as all remaining frontal slices. 
\end{defi}
\begin{defi}[orthogonal tensor \cite{KM11}]
	A tensor $\mathcal{Q}\in\mathbb{R}^{n\times n\times n_3}$ is orthogonal if $\mathcal{Q}^{\top}*\mathcal{Q}=\mathcal{Q}*\mathcal{Q}^{\top}=\mathcal{I}$.
\end{defi}

The t-product admits an analogue of the matrix singular value decomposition.
\begin{theorem}[t-SVD \cite{KM11}]\label{thm:tsvd}
	For any tensor $\mathcal{X}\in\mathbb{R}^{n_1\times n_2\times n_3}$, its tensor Singular Value Decomposition (t-SVD) is given by
	$\mathcal{X}=\mathcal{U}*\mathcal{S}*\mathcal{V}^{\top}$,
	where $\mathcal{U}\in\mathbb{R}^{n_1\times n_1\times n_3}$ and $\mathcal{V}\in\mathbb{R}^{n_2\times n_2\times n_3}$ are orthogonal tensors, and $\mathcal{S}\in\mathbb{R}^{n_1\times n_2\times n_3}$ is f-diagonal. If $\mathcal{X}$ has tensor tubal rank $r$ as defined in Definition~\ref{def:rank}, its skinny t-SVD is
	$\mathcal{X}=\mathcal{U}_r*\mathcal{S}_r*\mathcal{V}_r^{\top}$,
	where $\mathcal{U}_r$, $\mathcal{S}_r$, and $\mathcal{V}_r$ retain the first $r$ singular tubes.
\end{theorem}
\begin{defi}[tensor tubal rank \cite{KBHH13}]\label{def:rank}
	Given the t-SVD $\mathcal{X}=\mathcal{U}*\mathcal{S}*\mathcal{V}^{\top}$, the tensor tubal rank is the number of nonzero singular tubes in $\mathcal{S}$, i.e.,
	$\operatorname{rank}_t(\mathcal{X})=\#\{i:\mathcal{S}(i,i,:)\neq\bm{0}\}$.
\end{defi}

\begin{lem}\cite[Theorem 3.1]{Gra05}\label{Lem:CM}
	The eigenvalues of a circulant matrix $\mathscr{C}(\bm{z})$ generated by any vector $\bm{z}\in\mathbb{R}^n$ are given by the components of $F\bm{z}$, where $F$ is the DFT matrix.
\end{lem}
\begin{lem}\cite{BK22}\label{Lem:CC}
	For any vectors $\bm{x},\bm{y}\in\mathbb{R}^n$, one has $ \widehat{\bm{x}\bullet\bm{y}}=\hat{\bm{x}}\circ\hat{\bm{y}} $ and $ \|\bm{x}\|^2=\|\hat{\bm{x}}\|^2/n $, where $\circ$ denotes the Hadamard product.
\end{lem}



\section{MVCTD for Lightweight Online Adaptive Traffic Forecasting}\label{Sec:MVCTD}
This section develops the proposed Multi-View Coupled Tensor Decomposition (MVCTD) framework for online traffic prediction from imperfect multi-view streams. We first formalize the rolling forecasting task and then present the two stages of MVCTD. Stage I initializes a structured latent predictive representation from historical incomplete observations, while Stage II adapts this representation online as new partial traffic tensors arrive.

\subsection{Problem Formulation}\label{Sec:PD}
\begin{figure}[htbp]
	\centering
	\begin{tikzpicture}[
		scale=0.78,                
		>=Stealth, 
		thick,
		node distance=0.5cm,
		math font/.style={font=\normalsize\bfseries},       
		label font/.style={font=\scriptsize\bfseries\sffamily}, 
		tensor edge/.style={draw=gray!60, thick, line join=round},
		pixel/.style={draw=gray!30, very thin} 
		]
		
		\definecolor{colorInput}{RGB}{65, 105, 225}    
		\definecolor{colorTarget}{RGB}{255, 140, 0}    
		
		\newcommand{\drawTensor}[6]{
			\begin{scope}[shift={(#1)}]
				\pgfmathsetseed{\numexpr 100 + #6 \relax} 
				\def\cols{3} \def\rows{4} \def\sz{0.25} 
				\def\depthX{0.25} \def\depthY{0.15} 
				\def\w{\cols*\sz} \def\h{\rows*\sz}
				
				\filldraw[fill=#2!40, tensor edge] (\w, 0) -- ++(\depthX, \depthY) -- ++(0, \h) -- ++(-\depthX, -\depthY) -- cycle;
				\filldraw[fill=#2!20, tensor edge] (0, \h) -- ++(\w, 0) -- ++(\depthX, \depthY) -- ++(-\w, 0) -- cycle;
				
				\foreach \col in {0,...,\numexpr\cols-1} {
					\foreach \row in {0,...,\numexpr\rows-1} {
						\pgfmathparse{rnd}
						\ifdim \pgfmathresult pt > #5 pt 
						\filldraw[fill=#2!80, pixel] (\col*\sz, \row*\sz) rectangle ++(\sz, \sz);
						\else
						\filldraw[fill=white, pixel] (\col*\sz, \row*\sz) rectangle ++(\sz, \sz);
						\fi
					}
				}
				\draw[tensor edge] (0,0) rectangle (\w, \h);
				\node[below=0.15cm, math font, #2!80!black] at (\w/2, 0) {#3};
				
				\coordinate (#4_top) at (\w/2, \h); 
				\coordinate (#4_right) at (\w + \depthX, \h/2 + \depthY/2);
				\coordinate (#4_bot) at (\w/2, 0);
			\end{scope}
		}
		
		\begin{scope}[local bounding box=MainContent]
			
			\def\yOne{0}
			
			\drawTensor{0, \yOne}{colorInput}{\scriptsize$\mathcal{M}_1$}{t1_x1}{0.3}{1}   
			\node[font=\huge, colorInput!50] at (1.8, \yOne+0.5) {$\dots$};            
			\drawTensor{2.5, \yOne}{colorInput}{\scriptsize$\mathcal{M}_D$}{t1_xD}{0.2}{99}
			
			\draw[->, line width=1pt, gray!50] (t1_xD_right) -- ++(1.5, 0) node[midway, above, label font] {Forecast};  
			\drawTensor{5.0, \yOne}{colorTarget}{\scriptsize$\mathcal{X}_{D+1}$}{t1_pred}{0.0}{100}

			\def\yTwo{-2.8}
			
			\drawTensor{0, \yTwo}{colorInput}{\scriptsize $\mathcal{M}_1$}{t2_x1}{0.3}{1}
			\node[font=\huge, colorInput!50] at (1.8, \yTwo+0.5) {$\dots$};
			\drawTensor{2.5, \yTwo}{colorInput}{\scriptsize $\mathcal{M}_D$}{t2_xD}{0.2}{99}
			\drawTensor{5.0, \yTwo}{colorInput}{\scriptsize $\mathcal{M}_{D+1}$}{t2_xD1}{0.3}{100}
			
			\draw[->, line width=1pt, gray!50] (t2_xD1_right) -- ++(1.5, 0) node[midway, above, label font] {Forecast};
			\drawTensor{7.5, \yTwo}{colorTarget}{\scriptsize$\mathcal{X}_{D+2}$}{t2_pred}{0.0}{101}

			\def\yThree{-5.6}
			
			\drawTensor{0, \yThree}{colorInput}{\scriptsize $\mathcal{M}_1$}{t3_x1}{0.3}{1}
			\node[font=\huge, colorInput!50] at (1.8, \yThree+0.5) {$\dots$};
			\drawTensor{2.5, \yThree}{colorInput}{\scriptsize $\mathcal{M}_D$}{t3_xD}{0.2}{99}
			\drawTensor{5.0, \yThree}{colorInput}{\scriptsize $\mathcal{M}_{D+1}$}{t3_xD1}{0.3}{100}
			\drawTensor{7.5, \yThree}{colorInput}{\scriptsize $\mathcal{M}_{D+2}$}{t3_xD2}{0.3}{101}
			
			\draw[->, line width=1pt, gray!50] (t3_xD2_right) -- ++(1.5, 0) node[midway, above, label font] {Forecast};
			\drawTensor{10.0, \yThree}{colorTarget}{\scriptsize$\mathcal{X}_{D+3}$}{t3_pred}{0.0}{102}

			
			\draw[->, line width=1.2pt, colorInput!80] (t1_pred_top) ([yshift=1.2cm]t2_xD1_top) 
			-- node[midway, right, align=left, label font, xshift=2pt] {Data arrives\\(Incomplete)} (t2_xD1_top);
			
			\draw[->, line width=1.2pt, colorInput!80] (t2_pred_top) ([yshift=1.2cm]t3_xD2_top)
			-- node[midway, right, align=left, label font, xshift=2pt] {Data arrives\\(Incomplete)} (t3_xD2_top);
			
			\begin{scope}[on background layer]
				\draw[dashed, gray!30] (t1_x1_bot) -- (t2_x1_top);
				\draw[dashed, gray!30] (t2_x1_bot) -- (t3_x1_top);
				\draw[dashed, gray!30] (t1_xD_bot) -- (t2_xD_top);
				\draw[dashed, gray!30] (t2_xD_bot) -- (t3_xD_top);
				\draw[dashed, gray!20] (t2_xD1_bot) -- (t3_xD1_top);
			\end{scope}
			
		\end{scope}
		
		\coordinate (LegendPos) at ($(MainContent.south) + (0, -0.2)$);  
		
		\begin{scope}[shift={(LegendPos)}]
			\filldraw[fill=colorInput!80, pixel] (-4.0, -0.15) rectangle ++(0.3, 0.3);
			\node[right, font=\scriptsize\sffamily, text=black!70] at (-3.6, 0) {Observed data};  
			
			\filldraw[fill=white, pixel] (-1.2, -0.15) rectangle ++(0.3, 0.3);
			\node[right, font=\scriptsize\sffamily, text=black!70] at (-0.8, 0) {Missing data};
			
			\filldraw[fill=colorTarget!80, pixel] (1.8, -0.15) rectangle ++(0.3, 0.3);
			\node[right, font=\scriptsize\sffamily, text=black!70] at (2.2, 0) {Prediction};
		\end{scope}
		
	\end{tikzpicture}
	\caption{Schematic illustration of rolling online traffic forecasting with incomplete observations.}
	\label{fig:rolling_prediction}
\end{figure}
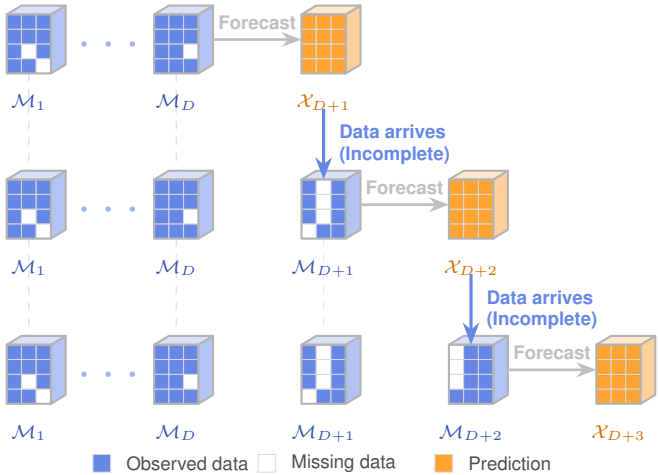

Let $\mathcal{X}_d \in \mathbb{R}^{N \times T \times V}$ denote the complete traffic state tensor on day $d$, where $N$, $T$, and $V$ are the numbers of road segments, daily time intervals, and traffic variables, respectively. The complete tensor is not fully observed in practice. Instead, the predictor receives $\mathcal{M}_d = \mathcal{P}_{\Omega_d}(\mathcal{X}_d)$, where $\Omega_d$ is the observed index set and $\mathcal{P}_{\Omega_d}$ denotes the projection that keeps entries in $\Omega_d$ and sets the remaining entries to zero. Given the partial history $\{\mathcal{M}_d\}_{d=1}^D$, the goal is to predict the next complete traffic tensor $\mathcal{X}_{D+1}$ and then repeat this procedure as new partial observations arrive, as illustrated in Fig. \ref{fig:rolling_prediction}.

This setting requires the model to use complementary information across traffic variables, recover predictive states from incomplete measurements, reduce the influence of non-recurrent anomalies, and update efficiently in a rolling deployment. Accordingly, MVCTD does not treat completion, anomaly suppression, and prediction as separate preprocessing and forecasting steps. Instead, it couples them within a single latent tensor formulation so that multi-view information can directly support forecasting under imperfect observations.

\subsection{Stage I: Offline Initialization}
Given the historical partial tensors $\{\mathcal{M}_d\}_{d=1}^D$ defined in Section \ref{Sec:PD}, the offline stage learns the initial components required for online prediction. It estimates regular latent traffic patterns, separates sparse abnormal disturbances, and initializes the temporal prediction model used in the subsequent streaming stage.

\textbf{Unified Optimization Framework.} Let $f_{\X}$ map a historical window of length $L$ to the current traffic state, i.e., $\hat{\mathcal{X}}_d = f_{\X}(\mathcal{X}_{d-L}, \ldots, \mathcal{X}_{d-1}; \Theta_{\X})$. Directly fitting this predictor to the raw stream is unreliable because the available observations are incomplete and may contain anomalous events. We therefore introduce a regular component $\mathcal{O}_d$ and an anomaly component $\mathcal{E}_d$, and learn the predictive state jointly with completion and anomaly separation:
\begin{equation}\label{eq:general_model}
	\begin{array}{cl}
		\min & \displaystyle \underbrace{\gamma \sum_{d=L+1}^{D} \| \mathcal{X}_d - f_{\X}(\mathcal{X}_{d-L}, \ldots, \mathcal{X}_{d-1}; \Theta_{\X}) \|^2}_{\text{Forecasting fidelity}} \\&+ \sum_{d=1}^{D} \big[  \mathcal{R}_1(\mathcal{O}_d) +  \mathcal{R}_2(\mathcal{E}_d)\big] \\
		\mbox{\rm s.t.} & \mathcal{X}_d = \mathcal{O}_d + \mathcal{E}_d, ~ \mathcal{P}_{\Omega_d}(\mathcal{X}_d) = \mathcal{M}_d, ~  d \in [D].
	\end{array}
\end{equation}
Here, $\mathcal{R}_1$ imposes structural regularization on the regular component $\mathcal{O}_d$, while $\mathcal{R}_2$ penalizes the anomaly component $\mathcal{E}_d$. The parameter $\gamma>0$ balances the forecasting fidelity term and the regularization terms.

\textbf{Coupled Tensor Factorization.} The regular component $\mathcal{O}_d$ should preserve the shared road-network structure across days while allowing different traffic variables to exhibit view-dependent temporal behavior. Since all daily traffic tensors are generated over the same road network, it is natural to use a common spatial dictionary to describe recurring spatial patterns. Meanwhile, day-to-day variations and view-specific temporal evolution are represented by day-specific latent temporal factors. The t-product further couples the factor tensors along the view dimension, enabling interactions among traffic variables to be encoded within the decomposition rather than treating each view independently. We therefore factorize the regular component as
\begin{equation*}
	\mathcal{O}_d = \mathcal{A} * \mathcal{Z}_d, \quad \forall d \in [D],
\end{equation*}
where $\mathcal{A} \in \mathbb{R}^{N \times R \times V}$ is a shared spatial dictionary and $\mathcal{Z}_d \in \mathbb{R}^{R \times T \times V}$ stores the day-specific latent temporal factors. This representation moves forecasting from the raw observation space to a structured latent space, where cross-view complementarities are encoded by the coupled dictionary and recurrent temporal patterns are represented more compactly.

The low-rank nature of $\mathcal{O}_d$ can be promoted through sparsity of the latent coefficients. Specifically, Yu et al. \cite[Theorem 3.2] {YDB26} established the equivalence between tensor tubal rank and group sparsity of the core coefficients under an orthogonal basis:
\begin{multline*}
	\operatorname{rank}_t(\mathcal{O}_d)=\min \left\{\|\mathcal{Z}_d\|_{\HS}^0: \mathcal{O}_d=\mathcal{A} * \mathcal{Z}_d, \mathcal{A}^{\top} * \mathcal{A}=\mathcal{I},\right. \\ \left. \mathcal{A} \in \mathbb{R}^{N \times R \times V}, \mathcal{Z}_d \in \mathbb{R}^{R \times T \times V}\right\},
\end{multline*}
where $\|\mathcal{Z}_d\|_{\HS}^0 = \sum_{r=1}^{R}\|\mathcal{Z}_d(r,:,:)\|^0$. 

Combining these physically motivated structures with theoretically grounded constraints, and using the $\ell_p$-norm relaxation for the $\ell_0$-norm, we reformulate problem \eqref{eq:general_model} as follows:
\begin{equation*}\label{eq:coupled_factorization}
	\begin{array}{cl}
		\min & \gamma \sum_{d=L+1}^{D} \| \mathcal{Z}_d - f_{\Z}(\mathcal{Z}_{d-L}, \ldots, \mathcal{Z}_{d-1}; \Theta_{\Z}) \|^2 \\
		&+ \sum_{d=1}^{D}\big[  \lambda_1\|\mathcal{Z}_d\|_{\HS}^p + \lambda_2\|\E_d\|_{\TF}^p\big]  \\
		\mbox{\rm s.t.} & \mathcal{X}_d = \mathcal{A} * \mathcal{Z}_d + \mathcal{E}_d, ~  \mathcal{P}_{\Omega_d}(\mathcal{X}_d) = \mathcal{M}_d,~d \in [D],\\
		& \mathcal{A}^{\top} * \mathcal{A}=\mathcal{I}.
	\end{array}
\end{equation*}
Here, $\lambda_1,\lambda_2>0$ are regularization parameters.
The first term measures prediction error in the latent space rather than in the high-dimensional traffic tensor space. Since $R \ll N$ and $\mathcal{E}_d$ absorbs sparse corruptions, latent-space forecasting reduces the dimension of the prediction model and makes the learned dynamics less sensitive to abnormal observations. The penalty $\|\mathcal{Z}_d\|_{\HS}^p = \sum_{r=1}^{R}\|\mathcal{Z}_d(r,:,:)\|^p$ promotes a compact regular component with low tubal rank, while $\|\mathcal{E}_d\|_{\TF}^p = \sum_{n=1}^{N}\sum_{t=1}^{T}\|\mathcal{E}_d(n,t,:)\|^p$ induces group sparsity across the view dimension \cite{YDC26}, allowing correlated abnormal responses to be separated from regular traffic dynamics.

\textbf{Subset Autoregression.}
Standard AR models typically rely on a continuous window of lags (e.g., $1, \dots, p$), which often leads to over-parameterization and the curse of dimensionality, particularly when capturing long-term dependencies. To address this, we characterize the temporal dynamics of $\Z_d$ using an SAR strategy. This approach is parsimonious: by explicitly selecting a sparse set of physically meaningful time lags $\mL:=\{h_l\}_{l=1}^{\bar L}$ (e.g., daily or weekly cycles) rather than consecutive ones, the model captures multi-scale periodicities with fewer parameters. Here, $\bar L=|\mL|$ denotes the number of selected lags, and $L:=\max\mL$ denotes the minimum history length required by the forecasting loss. Consequently, the evolution of the latent tensor is formulated as
\begin{equation*}\label{eq:SVAR}
	\Z_d \approx \sum_{l=1}^{\bar L} w^l \Z_{d-h_l},
\end{equation*}
where $w^l \in \R$ is the autoregressive coefficient associated with lag $h_l$.

\textbf{Laplacian Convolutional Priors.} 
Traffic states usually contain both global periodicity and local temporal continuity. Inspired by \cite{CCC24}, these two patterns can be characterized by a circulant matrix nuclear norm and a temporal convolution penalty. By imposing on the traffic state itself, the corresponding priors take the form
\begin{equation*}\label{eq:LC-X}
	\lambda_3\|\mathscr{C}\diamond\mathcal{X}_d\|_* + \lambda_4\|\bm{\ell} \boxtimes \mathcal{X}_d\|_{\RF},
\end{equation*}
where $\lambda_3,\lambda_4>0$ are regularization parameters, and $$\bm{\ell} \triangleq (2\tau, \underbrace{-1, \ldots,-1}_{\tau}, 0, \ldots, 0, \underbrace{-1, \ldots,-1}_{\tau})^{\top} \in \mathbb{R}^T$$ is the temporal convolution kernel. The term $\|\mathscr{C}\diamond \mathcal{X}_d\|_* = \sum_{n=1}^N\sum_{v=1}^V\|\mathscr{C}(\mathcal{X}_d(n,:,v))\|_*$ promotes global low-rank periodicity, and $\|\bm{\ell} \boxtimes \mathcal{X}_d\|_{\RF}= \sum_{n=1}^N\sum_{v=1}^V \|\bm{\ell} \bullet \mathcal{X}_d (n,:,v)\|$ promotes local smoothness.

In MVCTD, however, these temporal priors are imposed on $\mathcal{Z}_d$ rather than directly on $\mathcal{X}_d$. Since prediction is performed in the latent space, regularizing $\mathcal{Z}_d$ directly encourages extrapolatable recurrent patterns. Moreover, $\mathcal{X}_d$ may include sparse anomalies represented by $\mathcal{E}_d$; smoothing the raw state could therefore mix abnormal disturbances into regular dynamics.

The following bound supports this latent-space regularization strategy by showing that temporal regularity of the recovered traffic state can be controlled through the latent coefficients.

\begin{theorem}\label{Thm:bound}
	Let $\mathcal{A} \in \mathbb{R}^{N \times R \times V}$ and $\mathcal{Z} \in \mathbb{R}^{R \times T \times V}$, and let $\bm{\ell} \in \mathbb{R}^{T}$ be a convolution kernel. If $\mathcal{X} = \mathcal{A} * \mathcal{Z}$, then the following inequality holds.
	\begin{itemize}
		\item[(1)] $\|\mathscr{C}\diamond\mathcal{X}\|_* \le \zeta\|\mathscr{C}\diamond\mathcal{Z}\|_*$;
		\item[(2)] $\|\bm{\ell} \boxtimes \mathcal{X}\|_{\RF} \leq \zeta \|\bm{\ell} \boxtimes \mathcal{Z}\|_{\RF}$.
	\end{itemize}
	Here $\zeta = \max_{r} \sum_{n=1}^{N} \sum_{v=1}^{V} |\mathcal{A}(n, r, v)|$.
\end{theorem}

Motivated by this bound and the robustness considerations above, we reformulate the Laplacian convolutional priors in the latent space as
\begin{equation*}\label{eq:LC-Z}
	\lambda_3\|\mathscr{C}\diamond\mathcal{Z}_d\|_* + \lambda_4\|\bm{\ell} \boxtimes \mathcal{Z}_d\|_{\RF}.
\end{equation*}

\textbf{Integrated Model Formulation.} 
Combining the components above yields the concrete objective for the offline initialization stage:
\begin{equation}\label{model:off}
	\begin{array}{cl}
		\min\limits_{\Theta}&  \displaystyle \underbrace{\gamma \sum_{d=L+1}^{D}\big\|\Z_d-\sum_{l=1}^{\bar L}w^l\Z_{d-h_l} \big\|^2}_{\text{Forecasting}} + \displaystyle \sum_{d=1}^{D} \big[\underbrace{\lambda_1\|\mathcal{Z}_d\|_{\HS}^p}_{\text{Low-rank}} \\&+\underbrace{\lambda_2\|\E_d\|_{\TF}^p}_{\text{Anomaly}}+ \underbrace{\lambda_3\|\mathscr{C}\diamond\mathcal{Z}_d\|_*}_{\text{Global periodicity}} + \underbrace{\lambda_4\|\bm{\ell} \boxtimes \mathcal{Z}_d\|_{RF}}_{\text{Local smoothness}} \big] \\
		\mbox{\rm s.t.} & \mathcal{X}_d = \mathcal{A} * \mathcal{Z}_d + \mathcal{E}_d, ~  \mathcal{P}_{\Omega_d}(\mathcal{X}_d) = \mathcal{M}_d, ~ d \in [D],\\
		&\mathcal{A}^{\top}*\mathcal{A}=\mathcal{I},
	\end{array}
\end{equation}
where $\Theta = \{\{\mathcal{X}_d, \mathcal{Z}_d, \mathcal{E}_d\}_{d=1}^D, \{w^l\}_{l=1}^{\bar L}, \mathcal{A}\}$ denotes the set of all learnable parameters.

\subsection{Stage II: Lightweight Online Adaptive Prediction}

After offline initialization, MVCTD operates in a rolling manner. At each new day, the current model first produces a prior forecast and then corrects the latent state once the partial observation becomes available. Only the current latent tensor is optimized iteratively, while the global model components are updated separately.

\textbf{Forecasting.} 
Before observing day $D+1$, MVCTD predicts the latent tensor $\mathcal{Z}_{D+1}^+$ using the learned SAR relation. With the current spatial basis $\mathcal{A}_D$ and autoregressive coefficients $w_D^l$, the prior traffic forecast is
\begin{equation}\label{eq:pre}
	\mathcal{X}_{D+1}^+ = \mathcal{A}_{D} * \mathcal{Z}_{D+1}^{+}, ~ \text{where} ~ \mathcal{Z}_{D+1}^{+} = \sum_{l=1}^{\bar L}w^l_{D}  \mathcal{Z}_{D+1-h_l}^{\star}.
\end{equation}

\textbf{Online Inference and Model Adaptation.}
When the partial observation $\mathcal{M}_{D+1} = \mathcal{P}_{\Omega_{D+1}}(\mathcal{X}_{D+1})$ is received, MVCTD refines the current latent state and updates the model parameters through three decoupled subproblems.

\begin{itemize}
	\item Latent state inference: The current latent tensor $\mathcal{Z}_{D+1}^{\star}$ is estimated by solving
	\begin{equation}\label{model:on}
		\begin{array}{cl}
			\min\limits_{\mathcal{Z},\mathcal{E}} & \lambda_1\|\mathcal{Z}\|_{\HS}^p +\lambda_2\|\mathcal{E}\|_{\TF}^p + \lambda_3\|\mathscr{C}\diamond\mathcal{Z}\|_* \\
			&+ \lambda_4\|\bm{\ell} \boxtimes \mathcal{Z}\|_{RF} + \gamma\|\mathcal{Z}-\mathcal{Z}_{D+1}^{+}\|^2 \\
			\mbox{\rm s.t.} & \mathcal{P}_{\Omega_{D+1}}(\mathcal{A}_D * \mathcal{Z} + \mathcal{E}) = \mathcal{M}_{D+1}.
		\end{array}
	\end{equation}

	\item SAR coefficients update:
	The SAR coefficients are updated by minimizing the cumulative latent prediction error:
	\begin{equation}\label{SL2:W}
		\min_{w^l} \sum_{d=L+1}^{D+1} \big\|\mathcal{Z}_d^{\star} - \sum_{l=1}^{\bar L}w^l  \mathcal{Z}_{d-h_l}^{\star} \big\|^2.
	\end{equation}
	
	\item Spatial basis update: The shared spatial basis is then updated to incorporate the latest recovered state:
	\begin{equation}\label{SL2:A}
		\min_{\mathcal{A}} \sum_{d=1}^{D+1} \| \mathcal{A} * \mathcal{Z}_d^{\star} + \mathcal{E}_d^{\star} - \mathcal{X}_d^{\star} \|^2, \quad \text{s.t.} \quad \mathcal{A}^{\top} * \mathcal{A} = \mathcal{I}.
	\end{equation}
\end{itemize}

\textbf{Recursion.}  
The updated parameters $\mathcal{A}_{D+1}$ and $w_{D+1}^l$ are used to roll the model forward and generate the next forecast:
\begin{equation*}
	\mathcal{X}_{D+2}^+ = \mathcal{A}_{D+1} * \mathcal{Z}_{D+2}^{+}, ~ \text{where} ~ \mathcal{Z}_{D+2}^{+} = \sum_{l=1}^{\bar L}w^l_{D+1}  \mathcal{Z}_{D+2-h_l}^{\star}.
\end{equation*}

\begin{figure*}[!htbp]
	\centering
\begin{tikzpicture}[
	>=Stealth,
	thick,
	font=\small,
	node distance=0.8cm and 1.0cm,
	box/.style={
		rectangle, rounded corners=4pt, draw,
		align=center,
		drop shadow={opacity=0.15, shadow xshift=2pt, shadow yshift=-2pt},
		line width=0.8pt,
		inner sep=8pt
	},
	cData/.style={box, draw=blue!50!black, fill=blue!5},
	cModel/.style={box, draw=orange!60!black, fill=orange!5},
	cOptim/.style={box, draw=green!40!black, fill=green!5},
	cParam/.style={box, draw=gray!60!black, fill=gray!5},
	container/.style={
		draw=gray!30, dashed, fill=gray!2, rounded corners=8pt, inner sep=0.5cm
	},
	title/.style={
		font=\bfseries\normalsize, color=gray!90, anchor=north west
	}
	]
	
	
	\node[cData, text width=3.5cm] (S1_Data) {
		\textbf{Historical Batch Data}\\[0.1em]
		\scriptsize
		Incomplete tensors $\{\mM_d\}_{d=1}^D$
	};
	
	\node[cParam, right=5.5cm of S1_Data, text width=5.2cm] (S1_Params) {
		\textbf{Initial Predictive Model Parameters}\\[0.1em]
		\scriptsize
		$\A_D, \{w_D^l\}_{l=1}^{\bar L}, \{\Z_d^{\star}\}_{d = 1}^D$
	};
	
	\draw[->, line width=1.2pt] (S1_Data) -- node[midway, align=center, yshift=2pt] {
		\textbf{Solve Problem \eqref{model:off} via ADMM}\\
	} (S1_Params);
	
	\begin{scope}[on background layer]
		\node[container, fit=(S1_Data)(S1_Params), label={[title]north west:Stage I: Offline Initialization}] (Box1) {};
	\end{scope}

	
	
	\node[cModel, below=2.5cm of S1_Data, anchor=west, text width=3.8cm] (S2_Predict) {
		\textbf{1. Forecasting}\\[0.1em]
		\scriptsize 
		$
		\left\{
		\begin{aligned}
			\mathcal{Z}_{D+1}^{+} &= \sum\nolimits_{l} w^l_{D}  \mathcal{Z}_{D+1-h_l}^{\star} \\
			\mathcal{X}_{D+1}^+ &= \mathcal{A}_{D} * \mathcal{Z}_{D+1}^{+}
		\end{aligned}
		\right.
		$
	};
	
	\node[cData, below=0.6cm of S2_Predict, text width=3.8cm, 
	double copy shadow={shadow xshift=3pt, shadow yshift=-3pt, opacity=0.3}] (S2_Observe) {
		\textbf{2. Streaming Data Arrival}\\[0.1em]
		\scriptsize
		Incomplete tensor $\mM_{D+1}$
	};

	\node[cOptim, right=1.2cm of $(S2_Predict.south east)!0.5!(S2_Observe.north east)$, 
	text width=5.8cm] (S2_Merged) {
		
		\textbf{3. Refined Predictive Model Parameters}\\[0.1em]
		
		\scriptsize
		\textbf{(i) Latent Factor Inference}\\
		Solve problem \eqref{model:on} via ADMM $\Rightarrow$ $\Z_{D+1}^{\star}$\\[0.3em]
		
		\hdashrule{5cm}{0.5pt}{2pt 2pt} \\[0.3em]
		
		\textbf{(ii) Closed-form Parameter Optimization}\\
		Update SAR coefficients $\Rightarrow$ $\{w_{D+1}\}_{l=1}^{\bar L}$\\
		Update spatial basis $\Rightarrow$ $\A_{D+1}$
		
	};

	
	\draw[->] (S2_Predict.east) -- ++(0.5,0) |- 
	node[pos=0.75, above, font=\tiny, text=orange!60!black] {}
	node[pos=0.75, below, font=\tiny, text=orange!60!black] {}
	($(S2_Merged.west)+(0, 0.4)$);
	
	\draw[->] (S2_Observe.east) -- ++(0.5,0) |- 
	node[pos=0.75, below, font=\tiny, text=blue!60!black] {} 
	($(S2_Merged.west)+(0, -0.4)$);
	
	\begin{scope}[on background layer]
		\node[
		container,
		fit=(S2_Predict)(S2_Observe)(S2_Merged),
		inner xsep=22pt,
		inner ysep=22pt,
		label={[title]north west:Stage II: Lightweight Online Adaptive Prediction}
		] (Box2) {};
	\end{scope}

	
	\coordinate (Left_Turn_X) at ($(S2_Predict.west) - (1.2, 0)$);  
	\coordinate (Inner_Turn_X) at ($(S2_Predict.west) - (0.6, 0)$); 
	
	\coordinate (Predict_In_Top) at ($(S2_Predict.west) + (0, 0.2)$); 
	\coordinate (Predict_In_Bot) at ($(S2_Predict.west) + (0, -0.2)$); 

	\path (Box1.south) -- coordinate[midway] (mid_corridor) (Box2.north);
	
	\draw[->, line width=1.5pt, color=orange!80!black, rounded corners=8pt] 
	(S1_Params.south) 
	|- (Left_Turn_X |- mid_corridor)      
	|- (Predict_In_Top);

	\coordinate (LoopBottom) at ($(Box2.south) + (0, 0.25)$);
	
	\draw[->, dashed, line width=1.2pt, gray!80!black, rounded corners=8pt] 
	(S2_Merged.south) 
	|- (LoopBottom -| Inner_Turn_X)       
	|- (Predict_In_Bot);                  
	
	\node[font=\scriptsize\bfseries, text=black, below] 
	at ($(S2_Merged.south |- LoopBottom)!0.5!(Inner_Turn_X |- LoopBottom)$) 
	{Recursive Rolling to $D+1$};
	
\end{tikzpicture}
	\caption{The proposed MVCTD framework.}
\label{fig:framework}
\end{figure*}
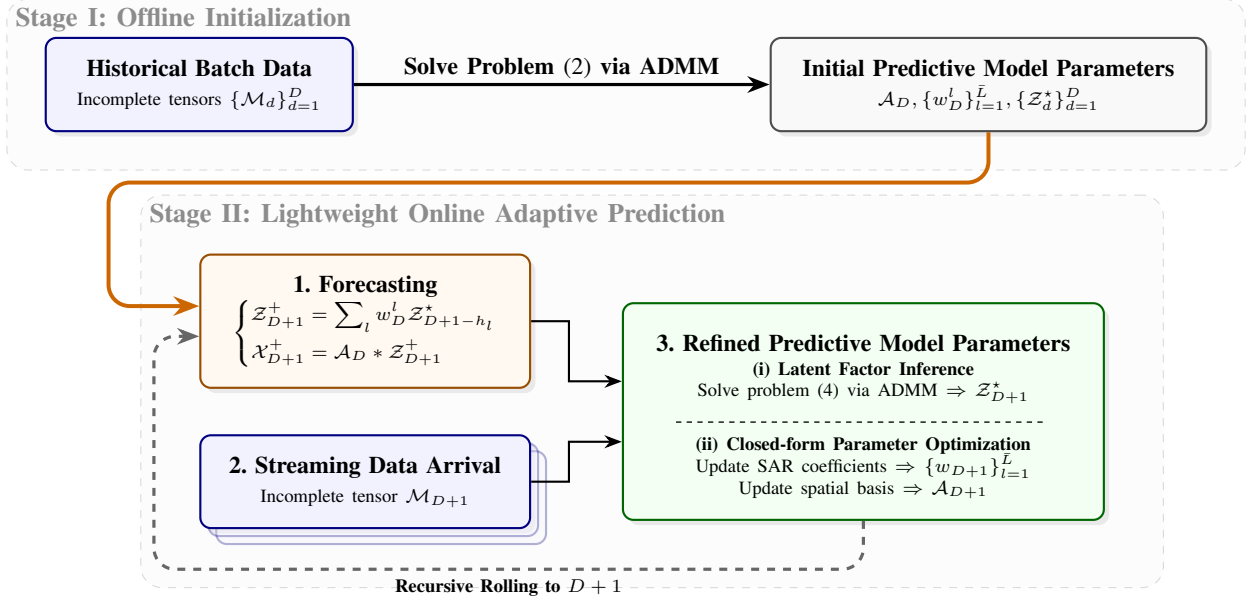

\begin{remark}[Computational efficiency]
	The efficiency of the online stage comes from separating current state inference from global parameter updates. When a new partial observation arrives, the previously inferred latent tensors $\{\mathcal{Z}_1^{\star}, \dots, \mathcal{Z}_D^{\star}\}$ are kept fixed, and iterative optimization is performed only for the current day latent inference problem. The SAR coefficients and the shared spatial basis are then updated through closed-form steps using summarized historical information. In this way, MVCTD avoids re-optimizing the full historical sequence during online adaptation.
\end{remark}

The full workflow of MVCTD, from offline initialization to online adaptive prediction, is summarized in Fig. \ref{fig:framework}.

\subsection{Theoretical Unification of Multi-View Coupled Tensor Decomposition}
To clarify the representational scope of the proposed t-product model, we relate it to standard matrix decompositions of individual traffic views. The following result shows that the coupled tensor representation can embed both independently decomposed views and the special case where all views share a common spatial basis.

\begin{theorem}\label{Thm:tproduct}
	Let $\mathcal{X} \in \mathbb{R}^{N \times T \times V}$ be a third-order tensor. 
	Suppose that for each view $v \in [V]$, the corresponding frontal slice admits a matrix decomposition
	\[
	X^{(v)}:=\X(:,:,v) = A_v Z_v,
	\]
	where $A_v \in \mathbb{R}^{N \times r_v}$ and $Z_v \in \mathbb{R}^{r_v \times T}$. 
	Then there exist an integer $p \le \sum_{v=1}^V r_v$ and two tensors 
	$\mathcal{A} \in \mathbb{R}^{N \times p \times V}$ and 
	$\mathcal{Z} \in \mathbb{R}^{p \times T \times V}$ such that
	\[
	\mathcal{X} = \mathcal{A} * \mathcal{Z}.
	\]
	
	More specifically, the following two cases can be embedded in this t-product form:
	\begin{itemize}
		\item Independent-view case. 
		When the views are allowed to have distinct spatial bases, one may set 
		$p = \sum_{v=1}^V r_v$ and construct the Fourier domain slices for each $ k \in [V] $ as
		\[
		\hat{A}^{(k)} = \big[ A_1 \; A_2 \; \cdots \; A_V \big], \quad
		\hat{Z}^{(k)} =
		\begin{bmatrix}
			\omega_k^{0} Z_1 \\
			\vdots \\
			\omega_k^{V-1} Z_V
		\end{bmatrix},
		\]
		where $\omega_k$ denotes the corresponding Fourier phase factor. 
		
		\item Shared-basis case.  
		If all views share an identical spatial basis, namely $A_v = A$ for all $v$, 
		the representation reduces to $p = r_1$. In this case, $\mathcal{A}$ reduces to the form
		\[
		A^{(1)} = A, \quad A^{(v)} = 0, \; v \neq 1,
		\]
		and the frontal slices of $\mathcal{Z}$ satisfy $Z^{(v)} = Z_v$. 
		Consequently, the coupled tensor model recovers the classical joint matrix decomposition form.
	\end{itemize}
\end{theorem}

This result shows that the proposed t-product representation covers both fully view-specific and shared-basis multi-view decompositions. Hence, MVCTD can model traffic variables that share common road-network structure while retaining view-dependent temporal patterns.

\section{Optimization Algorithm} \label{Sec:Alg}
To solve the non-convex offline model in \eqref{model:off}, we develop an Alternating Direction Method of Multipliers (ADMM) based alternating algorithm. Two auxiliary variables, $\Y_d = \Z_d$ and $\mC_d = \Z_d$, are introduced to decouple the latent temporal regularization and SAR fitting terms from the update of $\Z_d$. The resulting augmented Lagrangian function is given by
\begin{equation*}
	\begin{aligned}
		&\mathbb{L}_{\eta}\big(\X_d,\A,\Z_d,\Y_d,\mC_d,\E_d,w^l;\P_d,\Q_d,\mR_d\big)=\sum_{d=1}^{D}\big[ \lambda_1\|\Z_d\|_{\HS}^p \\
		&+\lambda_2\|\E_d\|_{\TF}^p + \lambda_3\|\mathscr{C}\diamond\Y_d\|_* + \lambda_4\|\bm{\ell} \boxtimes \Y_d\|_{RF}\big]  +  \\ & \gamma\sum_{d=L+1}^{D}\big\|\Z_d-\sum_{l=1}^{\bar L}w^l  \mC_{d-h_l} \big\|^2+\sum_{d=1}^{D}\big[\big\langle \A * \Z_d + \E_d - \X_d,\\& \P_d  \big\rangle + \frac{\eta_1}{2}\|\A * \Z_d + \E_d - \X_d\|^2\big] + \sum_{d=1}^{D}\big[\big\langle \Y_d - \Z_d, \Q_d  \big\rangle + \\ 
		&\frac{\eta_2}{2}\|\Y_d - \Z_d\|^2\big] + \sum_{d=1}^{D}\big[\big\langle \mC_d - \Z_d, \mR_d  \big\rangle + \frac{\eta_3}{2}\|\mC_d - \Z_d\|^2\big].
	\end{aligned}
\end{equation*}
Here, $\P_d$, $\Q_d$, and $\mR_d$ are dual variables, and $\eta_i>0,~i\in[3]$ are penalty parameters.
The primal variables are updated by minimizing this augmented Lagrangian with respect to one block at a time, followed by the dual variable updates. The main subproblems are derived below.

\textbf{Computing $\X_d$:} The sub-problem to update $\X_d$ is
\begin{equation}\label{CP:Xd}
	\mathop{\arg\min}\limits_{\P_{\Omega_d}(\X_d) = \mM_d}    \frac{\eta_1}{2}\|\A * \Z_d + \E_d - \X_d + \P_d/\eta_1\|^2.
\end{equation}
Clearly, the optimal solution to the $\X_d$ sub-problem is given by
\begin{equation}\label{SL:Xd}
	\X_d = \P_{\Omega_d^c}\big(\A * \Z_d + \E_d + \P_d/\eta_1\big) + \mM_d.
\end{equation}

\textbf{Computing $\A$:} The sub-problem to update $\A$ is
\begin{equation}\label{CP:A}
\begin{aligned}
	&\mathop{\arg\min}\limits_{\A^{\top}*\A = \I}~\sum_{d=1}^{D}\frac{\eta_1}{2}\|\A * \Z_d + \E_d - \X_d + \P_d/\eta_1\|^2\\
	=&\mathop{\arg\max}\limits_{\A^{\top}*\A = \I}~\big\langle \A, \sum_{d=1}^{D}\big(\eta_1\X_d-\eta_1\E_d-\P_d\big)*\Z_d^{\top}\big\rangle.
\end{aligned}
\end{equation}
By using \cite[Lemma 4.1]{YDB26}, the optimal solution of \eqref{CP:A} is given by
\begin{equation}\label{SL:A}
	\A = \U_\A*\V_\A^{\top},
\end{equation}
where $ \U_\A $ and $ \V_\A $ are obtained by the skinny t-SVD decomposition: $\sum_{d=1}^{D}\big(\eta_1\X_d-\eta_1\E_d-\P_d\big)*\Z_d^{\top} = \U_\A * \mS_\A * \V_\A^{\top} $.

\textbf{Computing $\Z_d$:} The sub-problem to update $\Z_d$ is
\begin{equation}\label{CP:Zd}
\begin{aligned}
	\mathop{\arg\min}\limits_{\Z_d}~&\lambda_1\|\Z_d\|_{\HS}^p + \gamma \mathbbm{1}_{(L,\infty)}(d) \big\|\Z_d-\sum_{l=1}^{\bar L}w^l  \mC_{d-h_l} \big\|^2 \\&+\big\langle \A * \Z_d + \E_d - \X_d, \P_d  \big\rangle + \frac{\eta_1}{2}\|\A * \Z_d + \E_d \\&- \X_d\|^2  + \big\langle \Y_d - \Z_d, \Q_d  \big\rangle + \frac{\eta_2}{2}\|\Y_d - \Z_d\|^2+ \\
	&\big\langle \mC_d - \Z_d, \mR_d  \big\rangle +  \frac{\eta_3}{2}\|\mC_d - \Z_d\|^2\\
	=\mathop{\arg\min}\limits_{\Z_d}~&\lambda_1\|\Z_d\|_{\HS}^p +\frac{\lambda_{\Z}}{2}\big\|\Z_d - \tilde{\Z} / \lambda_{\Z}\big\|^2,
\end{aligned}
\end{equation}
where $ \lambda_{\Z} = 2\gamma\mathbbm{1}_{(L,\infty)}(d) + \eta_1+\eta_2 + \eta_3 $ and $ \tilde{\Z} = 2\gamma \mathbbm{1}_{(L,\infty)}(d)\sum_{l=1}^{\bar L}w^l  \mC_{d-h_l} + \eta_1\A^{\top}*(\X_d-\E_d)-\A^{\top}*\P_d + \eta_2\Y_d + \Q_d + \eta_3\mC_d + \mR_d  $. By \cite[Appendix A]{PC21}, for each $r \in [R]$, the optimal solution of \eqref{CP:Zd} is given by
\begin{equation}\label{SL:Zd}
	\Z_d\left(r,:,:\right)= \begin{cases}\operatorname{Prox}_{\frac{\lambda_1}{\lambda_{\Z}}|\cdot|^p}\left(\frac{\tilde{z}}{\lambda_{\Z}} \right)\cdot\tilde{\Z}\left(r,:,:\right)/\tilde{z}, & \text{if}~\tilde{z} \neq 0; \\ 0, & \text{if}~\tilde{z}=0,\end{cases}
\end{equation}
where $\tilde{z} = \|\tilde{\Z}\left(r,:,:\right)\|$.

\textbf{Computing $\Y_d$:} The sub-problem to update $\Y_d$ is
\begin{multline}\label{CP:Yd}
	\mathop{\arg\min}\limits_{\Y_d}~\lambda_3\sum_{r=1}^{R}\sum_{v=1}^{V} \|\mathscr{C}\big(\Y_d(r,:,v)\big)\|_*+\lambda_4\sum_{r=1}^{R}\sum_{v=1}^{V}\\\|\bm{\ell} \bullet \Y_d(r,:,v)\|^2 + \big\langle \Y_d - \Z_d, \Q_d  \big\rangle + \frac{\eta_2}{2}\|\Y_d - \Z_d\|^2.
\end{multline}
According to Lemmas \ref{Lem:CM} and \ref{Lem:CC}, each component $\bm{y}:=\Y_d(r,:,v)$ of the tensor $ \Y_d $ is updated independently as follows:
\begin{equation}\label{CP:Ydrv}
	\begin{aligned}
		&\mathop{\arg\min}\limits_{\bm{y}}~\lambda_3 \|\mathscr{C}\big(\bm{y}\big)\|_*+\lambda_4\|\bm{\ell} \bullet \bm{y}\|^2 + \frac{\eta_2}{2}\|\bm{y} - \bm{h}\|^2\\
		=&\mathop{\arg\min}\limits_{\bm{y}}~\lambda_3 \|\hat{\bm{y}}\|_1+\frac{\lambda_4}{T}\|\hat{\bm{\ell}} \circ \hat{\bm{y}}\|^2 + \frac{\eta_2}{2T}\|\hat{\bm{y}} - \hat{\bm{h}}\|^2,
	\end{aligned}
\end{equation}
where $\bm{h} = \Z_d(r,:,v) - \Q_d(r,:,v)/\eta_2$. On each $\hat{y}_i$, the optimization problem is given by
\begin{equation}\label{SL:Ydrvi}
	\begin{aligned}
		&\mathop{\arg\min}\limits_{\hat{y}_i}~\lambda_3T \big| \hat{y}_i\big|+\frac{2\lambda_4\hat{\ell}_i^2+\eta_2}{2}\big(\hat{y}_i-\frac{\eta_2\hat{h}_i}{2\lambda_4\hat{\ell}_i^2+\eta_2}\big)^2\\
		=&\operatorname{sgn}\big(\frac{\eta_2\hat{h}_i}{2\lambda_4\hat{\ell}_i^2+\eta_2}\big) \cdot \max\big\lbrace \big| \frac{\eta_2\hat{h}_i}{2\lambda_4\hat{\ell}_i^2+\eta_2}\big|-\frac{\lambda_3T}{2\lambda_4\hat{\ell}_i^2+\eta_2},0\big\rbrace.
	\end{aligned}
\end{equation}

\textbf{Computing $\mC_d$:} The sub-problem to update $\mC_d$ is
\begin{multline}\label{CP:Cd}
	\mathop{\arg\min}\limits_{\mC_d}~\gamma \sum_{d=L+1}^{D}\big\|\Z_d-\sum_{l=1}^{\bar L}w^l  \mC_{d-h_l} \big\|^2+\big\langle \mC_d - \Z_d, \mR_d  \big\rangle\\ + \frac{\eta_3}{2}\|\mC_d - \Z_d\|^2.
\end{multline}
Let $J(\mC_d)$ denote the objective function in \eqref{CP:Cd}. The gradient $\nabla_{\mC_d} J(\mC_d) = 0$ is computed as
\begin{multline}
	\eta_3 (\mC_d - \Z_d) + \mR_d + \gamma \sum_{l=1}^{\bar L} \mathbbm{1}_{(L,D]}(d+h_l) \cdot \big[ -2w^l  \big( \Z_{d+h_l} - \\\sum_{j=1}^{\bar L}w^j  \mC_{d+h_l-h_j} \big) \big] = 0.
\end{multline}
This yields the update rule for $\mC_d$ given by
\begin{multline}\label{SL:Cd}
	\mC_d =  \big( \eta_3 \Z_d - \mR_d + 2\gamma \sum_{l=1}^{\bar L} \mathbbm{1}_{(L,D]}(d+h_l) w^l  \big( \Z_{d+h_l} - \\\sum_{j \ne l}w^j \mC_{d+h_l-h_j} \big) \big) / 
	 \big( \eta_3  + 2\gamma \sum_{l=1}^{\bar L} \mathbbm{1}_{(L,D]}(d+h_l) (w^l)^2 \big).
\end{multline}

\textbf{Computing $\E_d$:} The sub-problem to update $\E_d$ is
\begin{equation}\label{CP:Ed}
	\begin{aligned}
		&\mathop{\arg\min}\limits_{\E_d}~\lambda_2\|\E_d\|_{\TF}^p +\frac{\eta_1}{2}\|\A * \Z_d + \E_d - \X_d + \P_d/\eta_1\|^2.
	\end{aligned}
\end{equation}
Analogous to the update of $\Z_d$, for each $n \in [N],~t \in [T]$, the optimal solution of \eqref{CP:Ed} is given by
\begin{equation}\label{SL:Ed}
	\E_d\left(n,t,:\right)= \begin{cases}\operatorname{Prox}_{\frac{\lambda_2}{\eta_1}|\cdot|^p}\left(\tilde{e}\right)\cdot\tilde{\E}\left(n,t,:\right)/\tilde{e}, & \text{if}~\tilde{e} \neq 0; \\ 0, & \text{if}~\tilde{e}=0,\end{cases}
\end{equation}
where $ \tilde{\E} = \X_d - \A * \Z_d -\P_d/\eta_1 $ and $\tilde{e} = \|\tilde{\E}\left(n,t,:\right)\|$.

\textbf{Computing $w^l$:} The sub-problem to update $w^l$ is
\begin{equation}\label{CP:W}
	\mathop{\arg\min}\limits_{w^l}~\sum_{d=L+1}^{D}\big\|\Z_d-\sum_{l=1}^{\bar L}w^l  \mC_{d-h_l} \big\|^2.
\end{equation}
Let $\bm{w} = [w^1,\ldots,w^{\bar L}]^{\top} \in \mathbb{R}^{\bar L}$ denote the coefficient vector. This is a linear least squares problem with the closed-form solution given by
\begin{equation}\label{SL:W}
	\bm{w} = F^{-1}\bm{g},
\end{equation}
where the entries of $F \in \mathbb{R}^{\bar L \times \bar L}$ and  $\bm{g} \in \mathbb{R}^{\bar L}$ are computed by $F(k,l) = \sum_{d=L+1}^{D} \langle \mC_{d-h_k}, \mC_{d-h_l} \rangle$ and $\bm{g}(k) = \sum_{d=L+1}^{D} \langle \Z_d, \mC_{d-h_k} \rangle$, respectively.

Based on the derived update rules, we summarize the complete batch ADMM procedure for the offline initialization in Algorithm \ref{Alg:offline}.

\begin{algorithm}[!htbp]
	\caption{Stage I: Offline Initialization via Batch ADMM.}\label{Alg:offline}
	\KwIn{Incomplete tensors $\{\mathcal{M}_d\}_{d=1}^D$, regularization parameters $\gamma$ and $\{\lambda_i\}_{i=1}^4$. Set $k\leftarrow 0$.} 
	
	Update $\X_d^{k+1}$ via \eqref{SL:Xd}; \label{Alg:offline-step1}
		
	Update $\A^{k+1}$ via \eqref{SL:A};
		
	Update $\Z_d^{k+1}$ via \eqref{SL:Zd};
		
	Update $\Y_d^{k+1}$ via \eqref{SL:Ydrvi};
		
	Update $\mC_d^{k+1}$ via \eqref{SL:Cd};
		
	Update $\E_d^{k+1}$ via \eqref{SL:Ed};
		
	Update ${w^l}^{k+1}$ via \eqref{SL:W};
		
	Update Lagrangian multipliers and penalty parameters via
	\begin{equation}\label{SL:Lag}
		\left\{\begin{array}{l}
			\P_d^{k+1} = \P_d^k + \eta_1^k(\A^{k+1} * \Z_d^{k+1} + \E_d^{k+1} - \X_d^{k+1}); \\
			\Q_d^{k+1} = \Q_d^k + \eta_2^k(\Y_d^{k+1} - \Z_d^{k+1});\\
			\mR_d^{k+1} = \mR_d^k + \eta_3^k(\mC_d^{k+1} - \Z_d^{k+1});\\
			\eta_i^{k+1}=\beta \eta_i^k, i \in[3].
		\end{array}\right.
	\end{equation}

	If a termination criterion is met, set $ \X_d^{\star} := \X_d^{k+1} $, $\E_d^{\star} := \E_d^{k+1}$, $\A_D:=\A^{k+1}$, $ \Z_d^{\star}:= \Z_d^{k+1}$, $ w_D^l:= {w^l}^{k+1}$; else, set $k\leftarrow k+1$, return to Step \ref{Alg:offline-step1};
		
	\KwOut{Recover tensors $\{\X_d^{\star}\}_{d=1}^D$, anomaly tensor $\{\E_d^{\star}\}_{d=1}^D$, dictionary tensor $\A_D$, latent tensors $\{\Z_d^{\star}\}_{d=1}^D$, and autoregressive coefficients $\{w_D^l\}_{l=1}^{\bar L}$.}
\end{algorithm}

Although Algorithm \ref{Alg:offline} provides a reliable initialization, applying it repeatedly after each new observation would be too expensive for streaming prediction. We therefore use a lightweight online procedure. The prior $\mathcal{Z}_{D+1}^{+}$ from \eqref{eq:pre} is used as a warm start, and ADMM is run only for the current day inference problem \eqref{model:on}. After the current latent tensor is obtained, the SAR coefficients and the spatial dictionary are refreshed through the closed-form problems \eqref{SL2:W} and \eqref{SL2:A}. This procedure is summarized in Algorithm \ref{Alg:online}. In practice, a full retraining or an incremental/global learning strategy such as \cite{DSD16} can be triggered when the validation error becomes large.

\begin{algorithm}[!htbp]
	\caption{Stage II: Lightweight Online Adaptive Prediction}\label{Alg:online}
	\KwIn{Tensors $\A_D$, $\{\Z_d^{\star}\}_{d=1}^D$, and $\{w_D^l\}_{l=1}^{\bar L}$ from Algorithm \ref{Alg:offline}.} 

	Predict $\mathcal{Z}_{D+1}^+$ and $\mathcal{X}_{D+1}^+$ via \eqref{eq:pre};
	
	If the prediction horizon is reached, terminate; otherwise, proceed to the next step;
	
	Initialize $\mathcal{Z}_{D+1}^0 = \mathcal{Z}_{D+1}^+$ and $\X_{D+1}^0 = \P_{\Omega^c}(\mathcal{Z}_{D+1}^+) + \mM_{D+1}$ \tcp*[r]{Warm-start initialization}
	
	Update $\mathcal{Z}_{D+1}^{\star}, \mathcal{E}_{D+1}^{\star}$ by solving model \eqref{model:on} via ADMM \tcp*[r]{Optimize for current day only}
	
	Sequentially update $w_{D+1}^l,~l\in[\bar L]$ and $\mathcal{A}_{D+1}$ via \eqref{SL2:W} and \eqref{SL2:A}; 
	
	$D\leftarrow D+1$;

	\KwOut{Predict stream $\mathcal{X}_{D+1}^+, \mathcal{X}_{D+2}^+, \ldots $.}
\end{algorithm}

\textbf{Complexity Analysis.}
Let $K_{\rm off}$ and $K_{\rm on}$ denote the numbers of ADMM iterations in Algorithms \ref{Alg:offline} and \ref{Alg:online}, respectively. The offline stage is dominated by t-product operations, FFT based temporal shrinkage, and the SAR least-squares update. The offline initialization requires $\mathcal{O}\big(K_{\rm off}\big(DVRT(N+\log T+\bar L^2)+\bar L^3\big)\big)$.
For small $R$ and $\bar L$, this cost scales nearly linearly with the historical data size.

At each online step, ADMM updates only the latent tensor of the newly arrived day, after which the model parameters are refreshed in closed form. The per-step complexity is
\begin{equation*}
		\mathcal{O}\big(K_{\rm on}(VNRT+RVT\log T)+\bar L^2RVT+\bar L^3\big).
\end{equation*}
Compared with batch re-optimization, the online update avoids the multiplicative factor $D$ in the iterative part, making MVCTD suitable for streaming prediction.

We now establish the convergence of the proposed Algorithm \ref{Alg:offline} through the following theorem.
\begin{theorem}\label{Thm:CA}
	Let $\{\Theta^k\}_{k=1}^{\infty}$ denote the sequence generated by Algorithm \ref{Alg:offline}, where
	$\Theta^k = \{\X_d^k, \A^k, \Z_d^k, \Y_d^k, \mC_d^k, \E_d^k, {w^l}^k, \P_d^k, \Q_d^k, \mR_d^k\}$.
	If the dual sequences $\{\Q_d^k\}_{k=1}^{\infty}$ and $\{\mR_d^k\}_{k=1}^{\infty}$ are bounded, then every accumulation point of $\{\Theta^k\}_{k=1}^{\infty}$ satisfies the Karush-Kuhn-Tucker (KKT) conditions of problem \eqref{model:off}.
\end{theorem}

\section{Numerical Experiments}\label{Sec:exp}
This section evaluates MVCTD on real-world traffic prediction tasks with incomplete multi-view observations. All experiments are conducted in MATLAB 2025b on a workstation with an Intel Core i5-12500H CPU (2.50 GHz) and 16 GB RAM.

\subsection{Experimental Settings}

\subsubsection{Implementation Details}
For all competing methods, we use the parameter settings recommended in their original papers. For MVCTD, the regularization parameters are empirically set to $\lambda_1 = 2e\text{-}1$, $\lambda_2 = 5e\text{-}1$, $\lambda_3 = 1e\text{-}2$, $\lambda_4 = 2e\text{-}1$, and $\gamma = 5e\text{-}1$. The penalty growth factor is set to $\beta = 1.15$. The SAR lag set is fixed as $\mL=\{7\}$, so that $\bar L = |\mL| = 1$ and $L=\max\mL=7$. The optimization is terminated if the following criterion is satisfied for five consecutive iterations:
\begin{multline*}
	\text{Err}:=\max\left\{ 
	\|\A^{k+1} * \Z_d^{k+1} - \A^k * \Z_d^k\|_\infty, 
	\right. \\ 
	\left. \|\E_d^{k+1} - \E_d^k\|_\infty,~\|\X_d^{k+1} - \X_d^k\|_\infty 
	\right\} < 1e\text{-}3.
\end{multline*}
The maximum number of iterations is set to 200.

\subsubsection{Traffic Datasets}
We use two real-world traffic datasets from the Caltrans Performance Measurement System (PeMS)\footnote{\url{https://pems.dot.ca.gov/}}. Both datasets contain volume, occupancy, and speed measurements at a 5-minute resolution, forming a three-view traffic tensor. Occupancy values are rescaled from $[0, 1]$ to $[0, 100]$ before evaluation.

\begin{itemize}
	\item PeMS-D4. This dataset contains measurements from 307 loop detectors in the San Francisco Bay Area (District 4) from January 1 to February 28, 2018. It contains $59$ days and $16,992$ time intervals, yielding a tensor of size $307 \times 16,992 \times 3$.
	
	\item PeMS-D8. This dataset contains measurements from 170 sensors in the San Bernardino area (District 8) from July 1 to August 31, 2016. It contains $62$ days and $17,856$ time intervals, yielding a tensor of size $170 \times 17,856 \times 3$.
\end{itemize}


\subsubsection{Evaluation Metrics}

Prediction accuracy is evaluated by Mean Absolute Percentage Error (MAPE) and Root Mean Square Error (RMSE):
\begin{equation*}
	\begin{aligned}
		\text{MAPE} &= \frac{1}{B} \sum_{i,j,k} \big| \frac{\mathcal{M}_{ijk} - \mathcal{X}_{ijk}}{\mathcal{M}_{ijk}} \big| \times 100, \\
		\text{RMSE} &= \sqrt{\frac{1}{B} \sum_{i,j,k} \big(\mathcal{M}_{ijk} - \mathcal{X}_{ijk}\big)^2},
	\end{aligned}
\end{equation*}
where $B$ is the number of evaluated entries, and $\mathcal{M}_{ijk}$ and $\mathcal{X}_{ijk}$ denote the ground-truth and predicted values at index $(i,j,k)$, respectively.

\subsubsection{Baseline Models}
We compare MVCTD with six representative baselines covering low-rank factorization, Bayesian modeling, dynamical system based modeling, and deep learning.

\begin{itemize}
	\item \textbf{HTMF} \cite{CZC25}. A temporal matrix factorization method that incorporates Hankel matrices into a low-dimensional latent space. It jointly optimizes factor matrices and Hankel structures via alternating minimization, enabling robust modeling of sparse and nonstationary traffic series.
	
	\item \textbf{BTMF} \cite{CS21}. A Bayesian temporal factorization model that integrates low-rank decomposition with a VAR process in a probabilistic graphical model. It captures global and local temporal consistencies for robust prediction and uncertainty estimation under missing data.
	
	\item \textbf{circDMDsp} \cite{WS23}. An anti-circulant dynamic mode decomposition method with sparsity constraints. It models traffic evolution as a high-dimensional nonlinear dynamical system and extracts interpretable spatio-temporal coherent structures for reconstruction and prediction.
	
	\item \textbf{PatchTST} \cite{NHNSK23}. A Transformer based model for multivariate time-series forecasting with a channel-independent architecture and patch-level segmentation. It reduces computational complexity while preserving local temporal semantics.
	
	\item \textbf{PreSTGNet} \cite{ZLW26}. A spatiotemporal graph neural network with a two-stage pre-training and fine-tuning framework. It learns transferable representations from traffic networks and adapts them to the target forecasting task.
	
	\textbf{AmTGBiM} \cite{XMJ26}. A Mamba based model for long-term traffic prediction that combines an adaptive masked spatial Transformer with dynamic graph attention. It captures long-range temporal dependencies and dynamic spatial correlations among traffic sensors.
\end{itemize}

\subsection{Online Traffic Forecasting}
We evaluate online forecasting performance on PeMS-D4 and PeMS-D8. Following the rolling prediction setting in Fig.~\ref{fig:rolling_prediction}, each experiment uses a 15-day historical window for training and a 5-day horizon for online prediction. To emulate sensor failures, 80\% of the input entries are randomly masked. This setting evaluates whether a predictor can remain accurate and efficient when the incoming multi-view traffic stream is incomplete. Since circDMDsp, PatchTST, PreSTGNet and AmTGBiM cannot directly handle incomplete inputs, they are evaluated only on fully observed data and used as reference baselines.

Tables \ref{tab:Pre_0_combined} and \ref{tab:Pre_80} report the forecasting results under fully observed and severely incomplete settings, respectively. When no entries are missing, MVCTD consistently obtains the lowest MAPE and RMSE on both datasets and across all three traffic views. The improvement is evident against both factorization-based and deep forecasting baselines. For example, on PeMS-D8, MVCTD reduces the best competing MAPE from 21.65 to 11.20 for volume and from 25.77 to 13.19 for occupancy, while also achieving the smallest speed error.

The superiority of MVCTD is preserved under severe data sparsity. With 80\% missing entries, MVCTD still ranks first on every reported metric in Table \ref{tab:Pre_80}. Its performance degradation from the fully observed case is also limited: the volume MAPE increases by 1.06 on PeMS-D4 and 0.95 on PeMS-D8, and the occupancy MAPE increases by less than one percentage point on both datasets. By comparison, HTMF and BTMF exhibit much larger prediction errors under the same missing rate, highlighting the robustness of MVCTD in highly incomplete online forecasting scenarios.

In addition to accuracy, MVCTD is computationally efficient. Under the fully observed setting, it requires 101 seconds on PeMS-D4 and 49 seconds on PeMS-D8, which is lower than all competing methods. Under the 80\% missing rate, the running time further decreases to 72 seconds and 38 seconds, respectively. Compared with the fastest baseline in Table \ref{tab:Pre_80}, MVCTD is about 26 times faster on PeMS-D4 and more than 33 times faster on PeMS-D8. These results indicate that MVCTD provides an effective accuracy--efficiency balance for online traffic prediction.

\begin{table}[htbp]
	\centering
	\caption{Comparison of forecasting performance with full observations.}
	\label{tab:Pre_0_combined}
	\resizebox{\linewidth}{!}{
		\begin{tabular}{@{} c c *{3}{c} *{3}{c} c @{}}
			\toprule
			\multirow{2}{*}{Dataset} & \multirow{2}{*}{Method} & \multicolumn{3}{c}{MAPE} & \multicolumn{3}{c}{RMSE} & \multirow{2}{*}{Time (s)} \\
			\cmidrule(lr){3-5} \cmidrule(lr){6-8}
			& & Volume & Occupancy & Speed & Volume & Occupancy & Speed & \\
			\midrule
			\multirow{7}{*}{PeMS-D4} 
			& MVCTD & \textbf{14.42} & \textbf{17.91} & \textbf{5.24} & \textbf{36.18} & \textbf{2.46} & \textbf{4.83} & \textbf{101} \\
			& HTMF & 32.41 & 59.13 & 6.80 & 55.13 & 3.00 & 5.54 & 4384 \\
			& BTMF & 36.52 & 64.11 & 7.43 & 67.78 & 3.31 & 5.99 & 663 \\
			& circDMDsp & 36.44 & 50.33 & 8.10 & 62.70 & 3.22 & 7.50 & 425 \\
			& PatchTST & 50.34 & 59.21 & 8.37 & 68.50 & 3.45 & 7.01 & 2686 \\
			& PreSTGNet & 38.08 & 51.09 & 9.18 & 70.94 & 3.81 & 7.44 & 23160 \\
			& AmTGBiM & 38.68 & 53.66 & 9.54 & 68.00 & 3.58 & 7.82 & 1418 \\
			\midrule
			\multirow{7}{*}{PeMS-D8}
			& MVCTD & \textbf{11.20} & \textbf{13.19} & \textbf{3.89} & \textbf{36.99} & \textbf{2.25} & \textbf{4.20} & \textbf{49} \\
			& HTMF & 23.96 & 25.77 & 5.30 & 51.48 & 2.64 & 4.70 & 3227 \\
			& BTMF & 24.21 & 33.09 & 5.87 & 57.22 & 2.87 & 5.19 & 1757 \\
			& circDMDsp & 21.65 & 26.78 & 6.06 & 50.43 & 2.64 & 5.79 & 289 \\
			& PatchTST & 26.85 & 32.49 & 5.96 & 58.30 & 2.85 & 5.38 & 2127 \\
			& PreSTGNet & 36.97 & 35.11 & 6.66 & 68.99 & 3.35 & 5.82 & 8409 \\
			& AmTGBiM & 23.93 & 35.66 & 7.55 & 61.47 & 3.03 & 6.56 & 855 \\
			\bottomrule
	\end{tabular}}
\end{table}

\begin{table}[htbp]
	\centering
	\caption{Comparison of forecasting performance under 80\% missing rate.}
	\label{tab:Pre_80}
	\resizebox{\linewidth}{!}{
		\begin{tabular}{@{} c c *{3}{c} *{3}{c} c @{}}
			\toprule
			\multirow{2}{*}{Dataset} & \multirow{2}{*}{Method} & \multicolumn{3}{c}{MAPE} & \multicolumn{3}{c}{RMSE} & \multirow{2}{*}{Time (s)} \\
			\cmidrule(lr){3-5} \cmidrule(lr){6-8}
			& & Volume & Occupancy & Speed & Volume & Occupancy & Speed & \\
			\midrule
			\multirow{3}{*}{PeMS-D4} 
			& MVCTD & \textbf{15.48} & \textbf{18.71} & \textbf{5.81} & \textbf{38.60} & \textbf{2.55} & \textbf{5.16} & \textbf{72} \\
			& HTMF & 35.68 & 52.79 & 7.18 & 58.62 & 2.94 & 5.68 & 2501 \\
			& BTMF & 36.59 & 65.04 & 7.60 & 69.45 & 3.33 & 6.07 & 1871 \\
			\midrule
			\multirow{3}{*}{PeMS-D8} 
			& MVCTD & \textbf{12.15} & \textbf{13.81} & \textbf{4.49} & \textbf{38.87} & \textbf{2.29} & \textbf{4.50} & \textbf{38} \\
			& HTMF & 27.46 & 37.09 & 6.28 & 57.33 & 3.83 & 5.19 & 3629 \\
			& BTMF & 23.87 & 32.11 & 5.83 & 57.30 & 2.85 & 5.12 & 1257 \\
			\bottomrule
	\end{tabular}}
\end{table}

\begin{figure*}[htbp]
	\centering
	\includegraphics[width=1\linewidth]{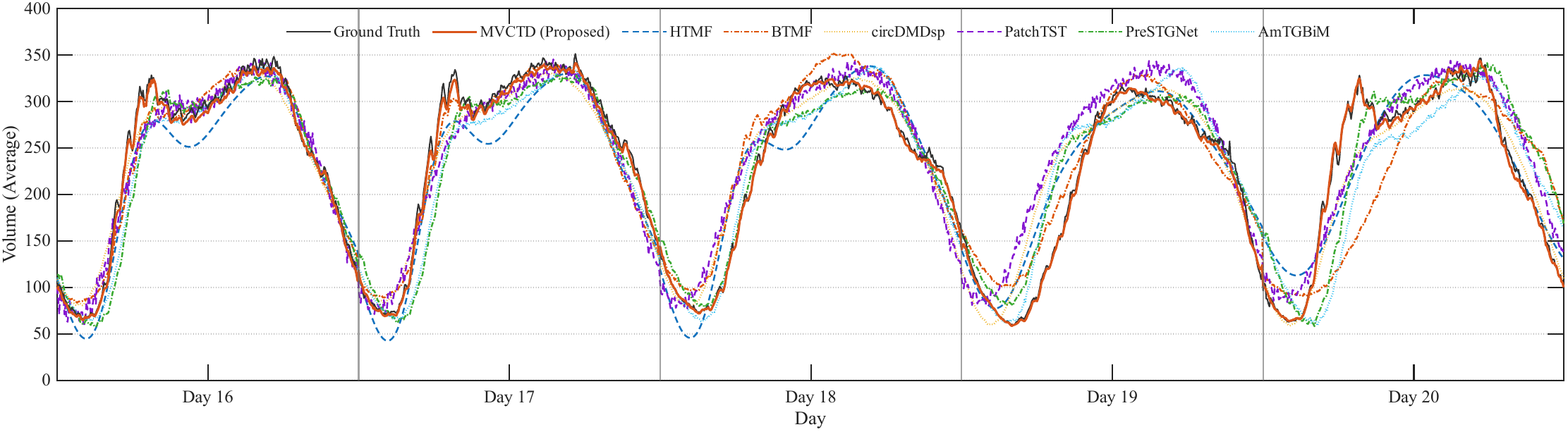}
	\caption{Time series prediction of average traffic volume on the PeMS-D8 dataset with full observations.}
	\label{fig:TimeSeries}
\end{figure*}

Fig. \ref{fig:TimeSeries} further compares the 5-day prediction trajectories of average traffic volume on PeMS-D8 with full observations. The ground-truth series presents a pronounced daily pattern, where sharp peak periods alternate with low volume intervals from Day 16 to Day 20. MVCTD reproduces this pattern with smaller amplitude and timing errors, particularly around the rapid rising and falling edges of the daily peaks. HTMF and BTMF capture the coarse periodic trend but deviate more noticeably near peak intervals, leading to either amplified fluctuations or attenuated peak values. The deep forecasting baselines, including PatchTST, PreSTGNet, and AmTGBiM, generally follow the global trend but produce smoother curves, which makes them less responsive to short term peak changes. circDMDsp shows a similar smoothing effect and misses part of the peak-to-valley variation. These trajectory-level differences are consistent with the volume MAPE and RMSE comparisons in Table \ref{tab:Pre_0_combined}.

\subsection{Model Analysis}

\subsubsection{Parameter Analysis}
\begin{figure*}[htbp]
	\centering
	\includegraphics[width=0.2\linewidth]{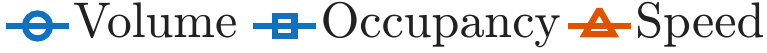}
	
	\includegraphics[width=0.2\linewidth]{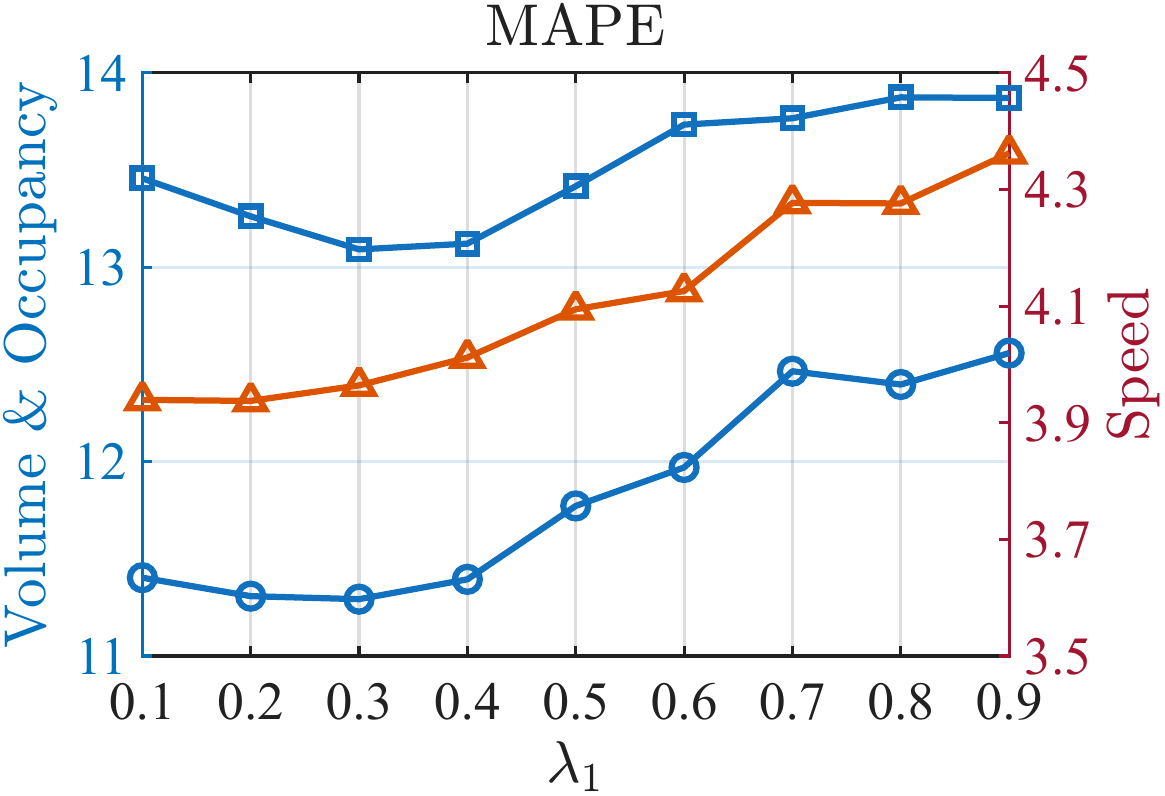}\hfill
	\includegraphics[width=0.2\linewidth]{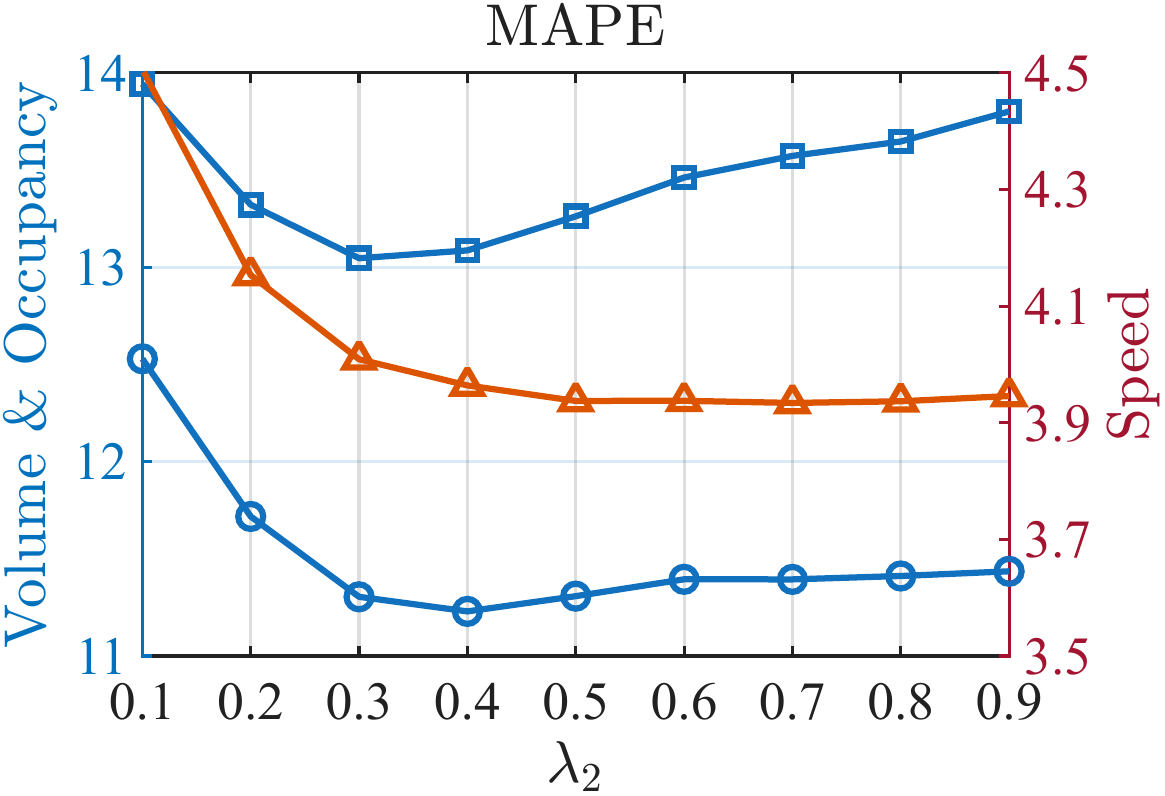}\hfill
	\includegraphics[width=0.2\linewidth]{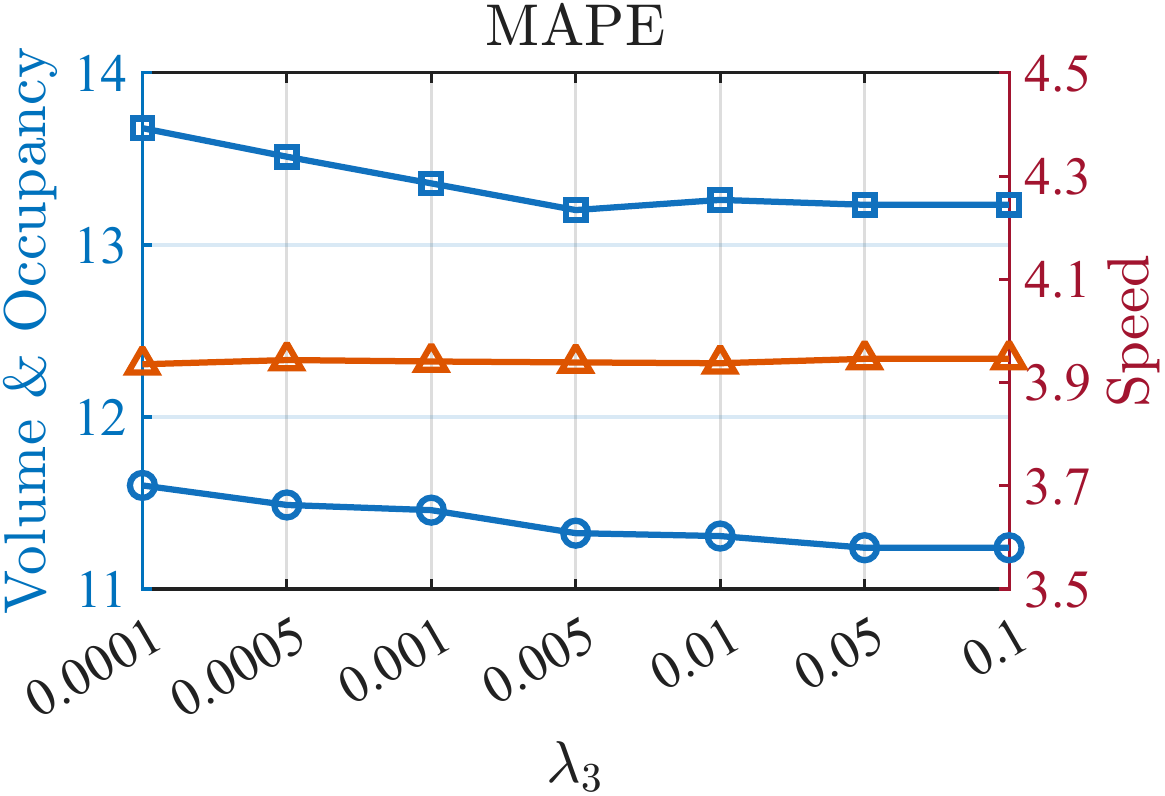}\hfill
	\includegraphics[width=0.2\linewidth]{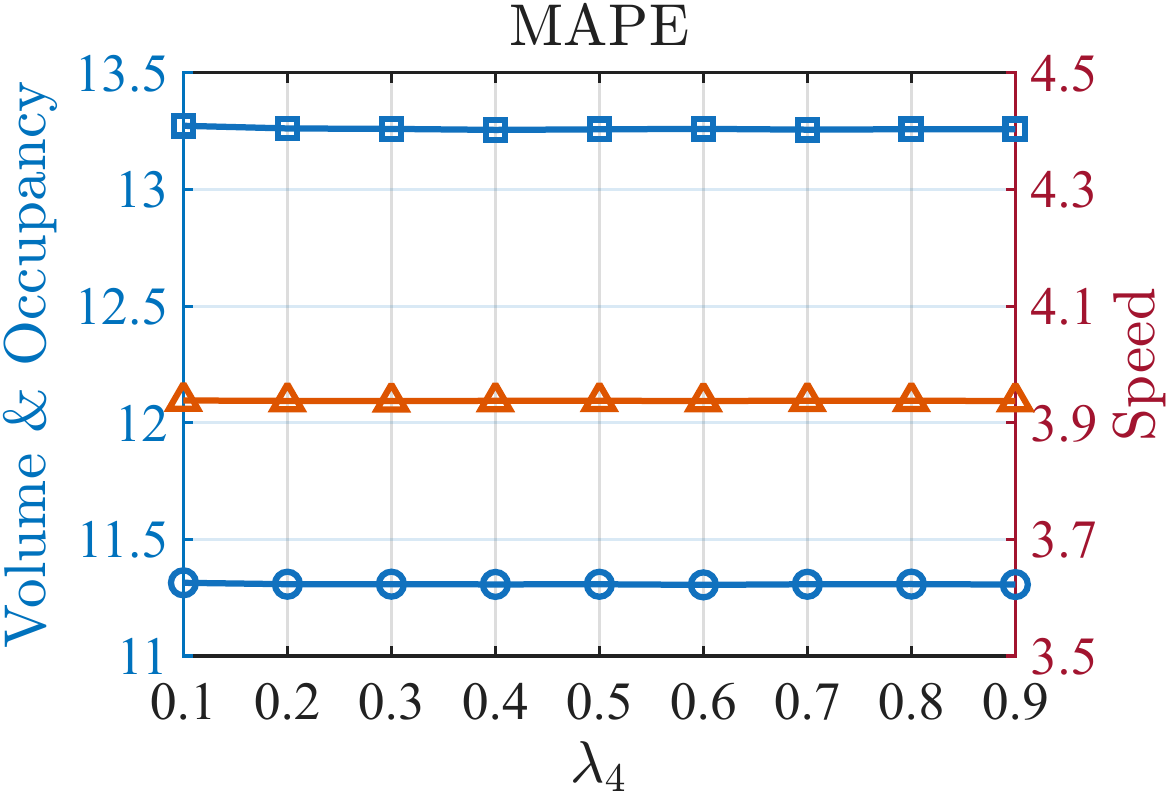}\hfill
	\includegraphics[width=0.2\linewidth]{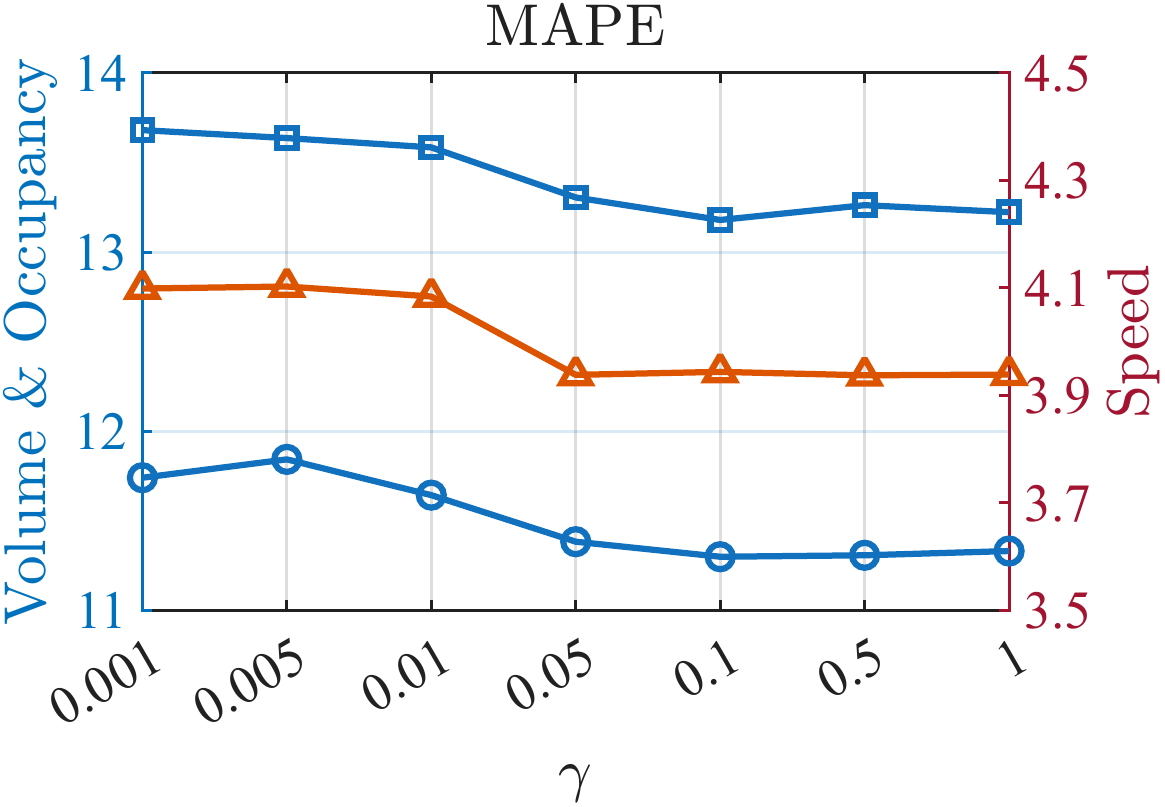}
	\caption{Parameter analysis of MVCTD with respect to $\lambda_i$ and $\gamma$ on the PeMS-D8 dataset under 40\% missing rate.}
	\label{fig:PA}
\end{figure*}

We evaluate the sensitivity of MVCTD to $\{\lambda_i\}_{i=1}^4$ and $\gamma$ on PeMS-D8 under a 40\% missing rate. As shown in Fig. \ref{fig:PA}, the model is more sensitive to $\lambda_1$ and $\lambda_2$, which control the latent low-rank and anomaly sparsity terms, respectively. The temporal regularization parameters $\lambda_3$ and $\lambda_4$ are relatively stable over a broad range, while $\gamma$ has a localized favorable region. Based on these results, we set $\lambda_1 = 2e\text{-}1$, $\lambda_2 = 5e\text{-}1$, $\lambda_3 = 1e\text{-}2$, $\lambda_4 = 2e\text{-}1$, and $\gamma = 5e\text{-}1$.

\subsubsection{Ablation Study}
Table \ref{tab:AS} reports the ablation results on PeMS-D8 with a 40\% missing rate. Removing the anomaly sparsity term ($\lambda_2$) leads to the largest degradation, increasing volume MAPE from 11.31 to 22.94. The global periodicity ($\lambda_3$) and local smoothness ($\lambda_4$) terms also contribute to the final performance, as removing either term weakens the prediction accuracy across the three views. The full model achieves the best or tied-best results on nearly all metrics, confirming the effectiveness of the proposed regularization design.

\begin{table}[htbp]
	\centering
	\caption{Ablation study of different regularization terms on the PeMS-D8 dataset under 40\% missing rate.}
	\label{tab:AS}
	\resizebox{\linewidth}{!}{
		\begin{tabular}{c c c c c c c}
			\toprule
			\multirow{2}{*}{Method} & \multicolumn{3}{c}{MAPE} & \multicolumn{3}{c}{RMSE} \\
			\cmidrule(lr){2-4} \cmidrule(lr){5-7}
			& Volume & Occupancy & Speed & Volume & Occupancy & Speed\\
			\midrule
			w/o Anomaly & 22.94 & 42.25 &  5.48 & 52.74 &  3.71 &  5.54    \\
			w/o Global Periodicity & 12.22 & 14.22 &  3.98 & 38.00 &  \textbf{2.25} &  4.29 \\
			w/o Local Smoothness & 11.61 & 14.07 &  3.97 & 37.85 &  2.27 &  4.28 \\
			Proposed & \textbf{11.31} & \textbf{13.26} &  \textbf{3.94} & \textbf{37.50} &  \textbf{2.25} &  \textbf{4.24}\\
			\bottomrule
	\end{tabular}}
\end{table}

\subsubsection{Effectiveness of Multi-View Modeling}
To evaluate the benefit of multi-view coupling, we compare MVCTD with a single-view variant trained independently for each traffic attribute. As shown in Table \ref{tab:MV}, MVCTD consistently outperforms the single-view variant on PeMS-D8, reducing speed MAPE from 4.67 to 4.49. On PeMS-D4, although the single-view variant obtains a slightly lower volume RMSE, MVCTD achieves lower MAPE values for all three views and reduces volume MAPE from 15.92 to 15.48. These results verify that coupling volume, occupancy, and speed provides useful complementary information under incomplete observations.

\begin{table}[htbp]
	\centering
	\caption{Performance comparison between single-view and multi-view modeling strategies.}
	\label{tab:MV}\resizebox{\linewidth}{!}{
	\begin{tabular}{cc ccc ccc}
		\toprule
		\multirow{2}{*}{Datasets} & \multirow{2}{*}{Strategy} & \multicolumn{3}{c}{MAPE (\%)} & \multicolumn{3}{c}{RMSE} \\
		\cmidrule(lr){3-5} \cmidrule(lr){6-8}
		&  & Volume & Occupancy & Speed & Volume & Occupancy & Speed \\
		\midrule
		\multirow{2}{*}{PeMS-D4}
		& Single-View & 15.92 & 19.87 & 6.02 & \textbf{37.86} & 2.61 & 5.37 \\
		& Proposed    & \textbf{15.48} & \textbf{18.71} & \textbf{5.81} & 38.60 & \textbf{2.55} & \textbf{5.16} \\
		\midrule
		\multirow{2}{*}{PeMS-D8}
		& Single-View & 12.24 & 14.03 & 4.67 & 39.17 & 2.40 & 4.75 \\
		& Proposed    & \textbf{12.15} & \textbf{13.81} & \textbf{4.49} & \textbf{38.87} & \textbf{2.29} & \textbf{4.50} \\
		\bottomrule
	\end{tabular}}
\end{table}

\subsubsection{Effectiveness of Subset Autoregression}

Table \ref{tab:SAR} compares different autoregressive lag structures. Using only the immediate previous day ($\{1\}$) gives the weakest performance, indicating that longer-term dependencies are important for traffic forecasting. The continuous window ($\{1, \dots, 7\}$) substantially reduces the errors, but the proposed sparse lag structure performs better on most metrics. Although its volume RMSE is slightly higher, it achieves the lowest errors in the other reported views. This result suggests that subset autoregression captures the dominant periodic dependence while avoiding redundant intermediate lags.

\begin{table}[htbp]
	\centering
	\caption{Impact of different autoregressive lag structures on prediction performance.}
	\label{tab:SAR}
	\resizebox{\linewidth}{!}{
		\begin{tabular}{c ccc ccc}
			\toprule
			\multirow{2}{*}{Time Lag Set} & \multicolumn{3}{c}{MAPE} & \multicolumn{3}{c}{RMSE} \\
			 \cmidrule(lr){2-4} \cmidrule(lr){5-7}
			  & Volume & Occupancy & Speed & Volume & Occupancy & Speed \\
			\midrule
			 $\{1\}$              & 22.74 & 26.77 &  5.62 & 60.62 &  3.03 &  5.44   \\
			 $\{1, \dots, 7\}$    &  12.32 & 14.56 &  4.18 & \textbf{37.15} &  2.28 &  4.33  \\
			 Proposed             & \textbf{11.31} & \textbf{13.26} &  \textbf{3.94} & 37.50 &  \textbf{2.25} &  \textbf{4.24}  \\
			\bottomrule
		\end{tabular}}
\end{table}

\subsubsection{Convergence Analysis}
\begin{figure}[htbp]
	\centering
	\begin{subfigure}[b]{0.333\linewidth}
		\centering
		\includegraphics[width=\linewidth]{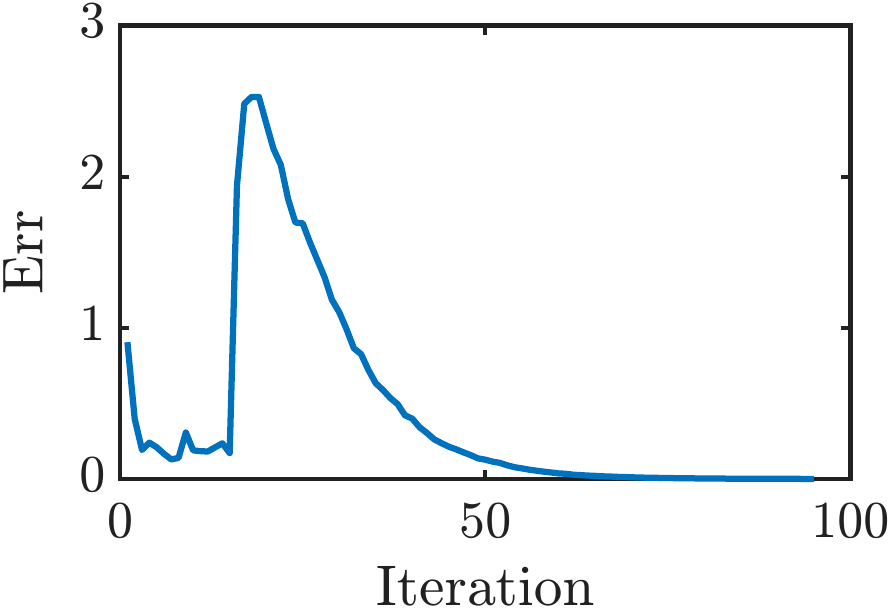}
		\caption{MR = $0\%$} 
	\end{subfigure}\hfill
	\begin{subfigure}[b]{0.333\linewidth}
		\centering
		\includegraphics[width=\linewidth]{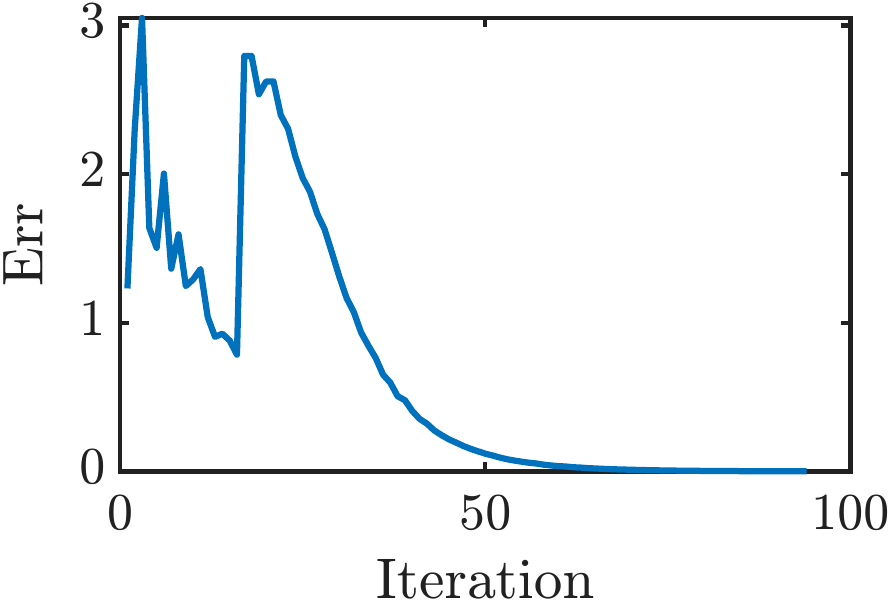}
		\caption{MR = $40\%$}
	\end{subfigure}\hfill
	\begin{subfigure}[b]{0.333\linewidth}
		\centering
		\includegraphics[width=\linewidth]{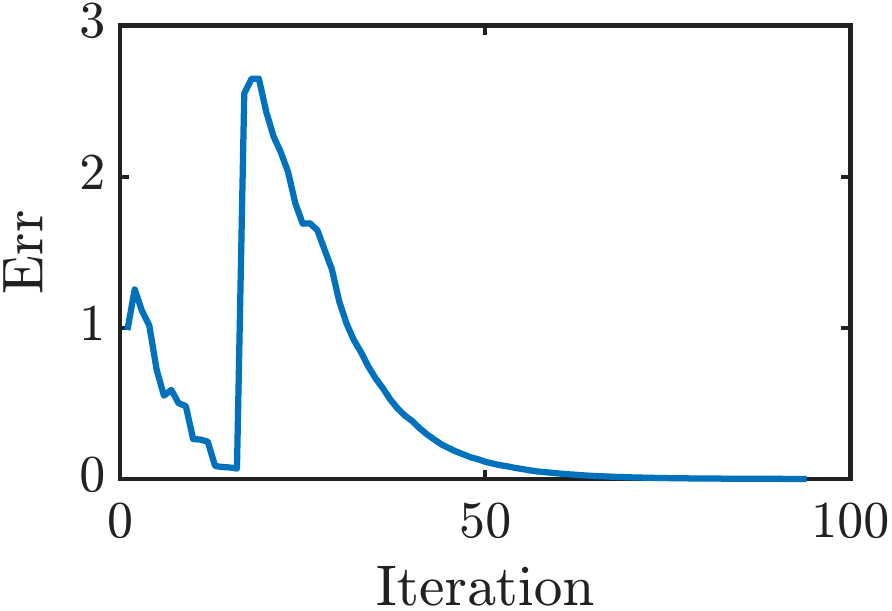}
		\caption{MR = $80\%$}
	\end{subfigure}
	\caption{Convergence curves of Algorithm \ref{Alg:offline} on the PeMS-D8 dataset under different Missing Rates (MR).}
	\label{fig:convergence_pems}
\end{figure}

Fig. \ref{fig:convergence_pems} shows the convergence curves of Algorithm \ref{Alg:offline} on PeMS-D8 under different missing rates. The error fluctuates during the first 20 iterations, which is expected for alternating minimization on a non-convex objective. It then decreases rapidly and reaches the stopping criterion within 90 iterations in all cases. This indicates that the proposed batch solver is practical for offline initialization.

\section{Conclusions}\label{Sec:con}

This paper proposed MVCTD, a multi-view coupled tensor decomposition model for online traffic prediction from incomplete and anomaly corrupted observations. MVCTD represents regular traffic states through a shared spatial dictionary and view-dependent latent temporal factors, so that complementary information among volume, occupancy, and speed can be exploited in a low-dimensional forecasting space. The model also integrates missing data recovery, anomaly separation, latent temporal regularization, and subset autoregression into a unified predictive formulation. To support streaming deployment, we developed a two-stage optimization scheme. The offline stage learns the initial coupled tensor representation from historical data, while the online stage refines only the current latent tensor and updates the remaining model parameters through lightweight closed-form steps based on summarized historical information. Experiments on PeMS-D4 and PeMS-D8 show that MVCTD achieves accurate and efficient prediction under severe missingness, and the ablation studies verify the contributions of anomaly modeling, temporal regularization, multi-view coupling, and subset autoregression. Future work will incorporate richer traffic-related factors, such as weather and network topology, and develop more adaptive online learning strategies for large-scale traffic systems.

\appendices
\ifCLASSOPTIONcaptionsoff
  \newpage
\fi



%
\bibliographystyle{IEEEtran}
\bibliography{References}

%




\end{document}